\documentclass[preprint,12pt]{elsarticle}

\usepackage[tbtags]{amsmath}
\usepackage{bm}
\usepackage{amsfonts}
\usepackage{booktabs}
\usepackage{adjustbox}
\usepackage{graphics}
\usepackage{graphicx}
\usepackage{psfrag}
\usepackage{eurosym}
\usepackage{pstricks}
\usepackage{array}
\usepackage{shortvrb}
\usepackage{epsf}
\usepackage{graphicx}
\usepackage{rotating}
\usepackage{float}
\usepackage{color}
\usepackage{tabularx}
\usepackage{multirow, makecell}
\usepackage{float}
\usepackage{colortbl}
\usepackage{ifpdf}
\usepackage{multicol}
\usepackage{todonotes}
\usepackage[ruled,vlined]{algorithm2e}

\DeclareMathOperator*{\argmin}{argmin}

\usepackage{hyperref}
\hypersetup{colorlinks=true, linkcolor=blue, anchorcolor=red, citecolor=blue, filecolor=red, urlcolor=red, pdfauthor=author}
            
\usepackage{subfigure}

\usepackage{geometry}
\usepackage{alphalph}
\usepackage{etoolbox}
\usepackage{pdfpages}

\usepackage[T1]{fontenc}

\usepackage{caption}

\newcommand{\comment}[1]{}

\makeatletter
\def\munderbar#1{\underline{\sbox\tw@{$#1$}\dp\tw@\z@\box\tw@}}
\makeatother

\patchcmd{\subequations}{\alph{equation}}{\alphalph{\value{equation}}}{}{}
\def \Noin {{\hskip 3pt \rm /\kern -9pt \in\hskip 1pt}}

\def \R {{\rm I\kern -2.2pt R\hskip 1pt}}

\allowdisplaybreaks

\begin{document}
\begin{frontmatter}

\title{A Quantile Neural Network Framework for Two-stage Stochastic Optimization}

\author[inst1]{Antonio Alc\'antara\corref{cor1}}
\ead{antalcan@est-econ.uc3m.es}
\author[inst1]{Carlos Ruiz}
\ead{caruizm@est-econ.uc3m.es}
\author[inst2]{Calvin Tsay}
\ead{c.tsay@imperial.ac.uk}
\cortext[cor1]{Corresponding author}

\affiliation[inst1]{organization={Department of Statistics, University Carlos III of Madrid}
            }

\affiliation[inst2]{organization={Department of Computing, Imperial College London}}

\begin{abstract}

Two-stage stochastic programming is a popular framework for optimization under uncertainty, where decision variables are split between first-stage decisions, and second-stage (or recourse) decisions, with the latter being adjusted after uncertainty is realized. These problems are often formulated using Sample Average Approximation (SAA), where uncertainty is modeled as a finite set of scenarios, resulting in a large ``monolithic'' problem, i.e., where the model is repeated for each scenario. The resulting models can be challenging to solve, and several problem-specific decomposition approaches have been proposed. An alternative approach is to approximate the expected second-stage objective value using a surrogate model, which can then be embedded in the first-stage problem to produce good heuristic solutions. In this work, we propose to instead model the \textit{distribution} of the second-stage objective, specifically using a quantile neural network. Embedding this distributional approximation enables capturing uncertainty and is not limited to expected-value optimization, e.g., the proposed approach enables optimization of the Conditional Value at Risk (CVaR). We discuss optimization formulations for embedding the quantile neural network and demonstrate the effectiveness of the proposed framework using several computational case studies including a set of mixed-integer optimization problems.
\end{abstract}

\begin{keyword}
Optimization under Uncertainty \sep \ Stochastic Programming \sep Neural Networks \sep Mixed-Integer Programming (MIP)
\end{keyword}

\end{frontmatter}

\section{Introduction}
\label{sec:intro}

Mathematical optimization provides a powerful framework for solving a wide range of decision-making problems, but conventional deterministic formulations rely on having exact estimates of involved model inputs and parameters. 
Therefore, when there is uncertainty associated with model inputs/parameters, stochastic programming (SP) approaches are preferred~\citep{birge2011introduction}. 
SP frameworks deal with solving optimization problems given \textit{known distributions} for uncertain inputs, rather than just their nominal values. 
As a result, SP has been applied in a variety of real-world optimization problems, such as in process systems engineering~\citep{li2021review}, supply chain optimization~\citep{govindan2018advances,santoso2005stochastic}, production scheduling~\citep{korpeouglu2011multi}, and unit commitment~\citep{shiina2004stochastic,van2018large,zheng2014stochastic}. 
We note that there are various other paradigms for optimization under uncertainty~\citep{powell2019unified}, such as chance-constrained optimization and robust optimization. 

Within SP, two-stage stochastic programming is a popular framework for handling exogenous (i.e., independent of the decision variables) uncertainties that are eventually realized. 
Specifically, two-stage stochastic programming divides decision variables into first-stage (``here-and now'') variables that must be decided before the realization of uncertainty, and second-stage (``wait-and-see'') ones that can be adjusted after the uncertainty is realized. 
This two-stage approach is popular owing to its flexibility, as it permits handling a wide range of uncertainties and allows adjustments to be made in the second stage once more information becomes available. 
Nevertheless, the primary challenge with two-stage stochastic programming relates to the computational expense in solving the resulting problems. 
In particular, uncertain parameters are represented using a set of scenarios, with each scenario corresponding to a possible realization of the parameters. The first-stage decisions are fixed across scenarios, and the expected value of the objective function can then be approximated by averaging over the full set of scenarios, known as the ``Sample Average Approximation'' (SAA). 
By choosing a risk metric such as Conditional Value-at-Risk (CVaR), two-stage stochastic programming can be extended to consider the second-stage risk as a function of the first-stage decision variables~\citep{rockafellar2000optimization,rockafellar2002conditional}. 

The SAA reformulation results in a deterministic, monolithic problem formulation that can therefore include repeated elements over a large number of scenarios, as well as linking constraints to enforce consistency of the first-stage decisions.
These problems can quickly become practically intractable if the underlying repeated model is large, the set of uncertain parameters is high-dimensional, and/or the number of scenarios considered is large (many scenarios are often needed to improve solution accuracy).
Various problem-specific approaches have been proposed to overcome these computational difficulties, as reviewed by~\citet{torres2022review}, such as decomposition approaches~\citep{ruszczynski1997decomposition}, dynamic optimization-based reformulation~\citep{tsay2017dynamic}, and scenario reduction strategies~\citep{heitsch2003scenario,mahmutougullari2018bounds}. 
In the first group, Benders decomposition, or the ``L-shaped'' method, has shown particular effectiveness when the second-stage problem is linear~\citep{van1969shaped}. 

Recently \citet{patel2022neur2sp} introduce \texttt{Neur2SP}, a framework for solving two-stage stochastic programs by using a surrogate model to approximate the second-stage problem. 
Specifically, the framework employs neural networks with ReLU activation functions to learn the \textit{expected} second-stage objective value as a function of the recourse variables for a set of scenarios. 
\citet{kronqvist2023alternating} later introduce an adaptive sampling technique for improving the accuracy of the surrogate model. 
Using the fact that ReLU neural networks can be encoded in mixed-integer optimization problems~\citep{huchette2023deep}, the learned second-stage surrogate model can be embedded in the overall optimization problem, replacing the SAA as an alternative approximation. 
In other words, the two-stage program is reformulated as the original first-stage problem where a neural network approximates the second stage.
This is a part of the growing literature on using embedded neural networks to replace intractable/unknown components of mixed-integer optimization, such as in constraint learning~\citep{fajemisin2023optimization} or distributional constraint learning~\citep{alcantara2023neural}.  
This \texttt{Neur2SP} framework was developed to specifically approximate the second-stage expected value, being therefore limited to risk-neutral decisions. Furthermore, the methodology can suffer from scalability of the data generation process or model size, as discussed later. 

In this work, we propose using quantile neural networks (QNNs) as the second-stage surrogate model. QNNs are effectively multi-output neural networks that can be similarly reformulated and embedded in optimization problems, and they enable quantifying the distributional aspect of the second-stage decisions. Specifically, the QNN approach enables learning the \textit{distribution} of the second-stage objective value as a function of the recourse variables, rather than just the \textit{expected value}. Two QNN model structures are presented: a standard unconstrained feed-forward NN, and an output-constrained NN that ensures the non-decreasing property of conditional quantile estimation.
We show that our framework is computationally efficient, with a fast data generation procedure and network training. Furthermore, once the QNN is embedded in the model, heuristic general solutions to the original problem can be obtained quickly ($<$1s in many cases) regardless of the size of the scenario set, as it is not included in the surrogate model. Regarding the quality of the solutions, the gap between the SAA and the proposed QNN approximation is small in many of our presented case studies, and the QNN obtains better results for large scenario sets when working under a time limit.

The remainer of the paper is organized as follows. Section~\ref{sec:background} presents background material on two-stage stochastic optimization, quantile neural networks, and mixed-integer formulations for embedding trained neural networks in optimization problems. 
In Section~\ref{sec:methodology} we then describe the proposed QNN-based stochastic programming framework, including procedures for the training and embedding of the QNN. 
Section~\ref{sec:cases_estudies} presents computational results for optimization of the expected value and risk-informed metrics (using CVaR) for several benchmark problems. 
We conclude in Section~\ref{sec:conclu}. 

\section{Background}
\label{sec:background}

\subsection{Two-stage Stochastic Optimization}
\label{sec:back_2sp}
A general representation of a two-stage stochastic problem is 
\begin{equation}\label{two_stage}
	\min_{\mathbf{X}\in \mathcal{X}} \ \mathbb{E}_{\boldsymbol{\xi}}[F(\mathbf{X},\xi)] = \min_{\mathbf{X}\in \mathcal{X}} \ c^{T}\mathbf{X} + \mathbb{E}_{\boldsymbol{\xi}}[V(\mathbf{X},\boldsymbol{\xi})]
\end{equation}
where $c\in\mathbb{R}^n$ and $\mathbf{X}\in \mathbb{R}^n$ define, respectively, the first-stage objective cost vector and the set of first-stage decision variables with feasible set $\mathcal{X}$ \citep{birge2011introduction}. $\boldsymbol{\xi}$ is the set of random parameters following a probability distribution $\mathcal{P}$ with support $\Psi$. For convenience, we are assuming that function $F(\mathbf{X},\boldsymbol{\xi})$ can be separated into a linear deterministic term $c^T\mathbf{X}$ and an arbitrary function $V(\mathbf{X}, \boldsymbol{\xi})$, as it is common to many applications, although this can be relaxed by setting $c=0$. Indeed, $V(\mathbf{X}, \boldsymbol{\xi})$ represents the second-stage value function $V:\mathcal{X}\times \Psi \rightarrow \mathbb{R}$, such that
\begin{equation}\label{second_Stage}
	V(\mathbf{X}, \boldsymbol{\xi}) = \min_{\mathbf{Y}\in\mathcal{Y}(\mathbf{X},\boldsymbol{\xi})} \ f(\mathbf{Y},\mathbf{X},\boldsymbol{\xi})
\end{equation}
where vector $\mathbf{Y} \in \mathbb{R}^m$ includes the second-stage recourse decision variables. Note that the first-stage decisions $\mathbf{X}$ and the random vector $\boldsymbol{\xi}$ parameterize the second stage objective function $f(\mathbf{Y},\mathbf{X},\boldsymbol{\xi})$ and feasibility region $\mathcal{Y}(\mathbf{X},\boldsymbol{\xi})$.

In general, problem (\ref{two_stage}) cannot be directly addressed unless important assumptions are made (e.g., linearity, independence, and normality). Hence, a common approach to tackle this problem is through its Sample Average Approximation (SAA). This is based on sampling the probability distribution $\mathcal{P}$ into a finite set of plausible scenarios $\boldsymbol{\xi}_\omega$ with $\omega=1,\dots,\Omega$, and by making the second-stage variables scenario dependent ($\mathbf{Y}_{\omega}$):
\begin{subequations}\label{eq:SAA_two_stage}
\begin{align}
&\min_{\mathbf{X},\mathbf{Y}_{\omega}} \ c^T \mathbf{X}+ \sum_{\omega=1}^{\Omega} \pi_{\omega} f(\mathbf{Y_\omega},\mathbf{X},\boldsymbol{\xi}_\omega)\\
&\rm{s.t.}\notag\\
& \quad \mathbf{X}\in \mathcal{X} \\
& \quad \mathbf{Y}_{\omega}\in\mathcal{Y}(\mathbf{X},\boldsymbol{\xi}_{\omega}), \quad \forall\ \omega=1,\dots,\Omega   
\end{align}
\end{subequations}
where $\pi_\omega$ is the probability associated with scenario $\omega$.

The larger the number of sampled scenarios $\Omega$, the better the underlying distribution $\mathcal{P}$ is characterized. However, problem (\ref{eq:SAA_two_stage}) involves $\omega$ duplicates of the second stage objective $f(\mathbf{Y},\mathbf{X},\boldsymbol{\xi})$ and feasibility region $\mathcal{Y}(\mathbf{X},\boldsymbol{\xi})$ and can easily become computationally intractable for nonlinear or mixed-integer linear formulations if this number is sufficiently large. 
Therefore, there is a natural tradeoff between tractability and accuracy of the SAA approach. 

Moreover, given $\mathbf{X}$, the objective function $F(\mathbf{X},\boldsymbol{\xi})$ can be viewed as a random variable. In this regard, there are many relevant applications where it is mandatory to account not only for its expected value but also for other aspects of its probability distribution (e.g., its variance, a given quantile, the worst case). This is in general addressed by extending problem (\ref{two_stage}) with mean-risk formulations \citep{ahmed2006convexity}, which combine the expected value and a risk measure $\mathcal{R}_{\boldsymbol{\xi}}$, resulting in:
\begin{equation}\label{eq:two_stage_risk}
	\min_{\mathbf{X}\in \mathcal{X}} \ \mathbb{E}_{\boldsymbol{\xi}}[F(\mathbf{X},\boldsymbol{\xi})]+\lambda \mathcal{R}_{\boldsymbol{\xi}}[F(\mathbf{X},\boldsymbol{\xi})]
\end{equation}
where $\mathcal{R}_{\boldsymbol{\xi}}:\mathcal{Z}\rightarrow\mathbb{R}$ and $\mathcal{Z}$ represents the space of all real random objective function values with $\mathbb{E}[F(\boldsymbol{\xi})]<\infty$, and $\lambda$ is a non-negative scalar that adjusts the trade-off between expectation and risk.

There exist several alternatives for $\mathcal{R}_{\boldsymbol{\xi}}$, although from a modeling perspective, the so-called \textit{coherent risk measures} are very interesting given to their convexity \citep{artzner1999coherent}. In particular, the conditional value-at-risk (CVaR), is one of the most popular risk measures used in practice. For a random variable $Z$, and a confidence level $\alpha \in (0,1)$, it can be defined as 
\begin{equation}\label{CVaR}
\text{CVaR}_\alpha = \mathbb{E}[Z|Z\geq\text{VaR}_\alpha(Z)]
\end{equation}  
where $\text{VaR}_\alpha(Z)$ is the value-at-risk for confidence level $\alpha$, i.e., the lower $\alpha$-quantile of the random variable $Z$.

Considering (\ref{eq:two_stage_risk}), using the CVaR as a risk measure is especially convenient for those cases where $\mathcal{X}$ is a convex set, $f(\mathbf{Y},\mathbf{X},\boldsymbol{\xi})$ is convex with respect to $\mathbf{X}$ and $\mathcal{Y}(\mathbf{X},\boldsymbol{\xi})$ is a polyhedral set, as the resulting model can be also approximated by a SAA-based approach \citep{rockafellar2000optimization,rockafellar2002conditional}:
\begin{subequations}\label{eq:SAA_two_stage_CVaR}
\begin{align}
&\min_{\mathbf{X},\mathbf{Y}_{\omega}} \ (1+\lambda)c^T \mathbf{X}+ \sum_{\omega=1}^{\Omega} \pi_{\omega} f(\mathbf{Y_\omega},\mathbf{X},\boldsymbol{\xi}_\omega)+\lambda \left(\nu+\frac{1}{1-\alpha}\sum_{\omega=1}^{\Omega}\pi_{\omega}\eta_{\omega}\right)\\
&\rm{s.t.}\notag\\
& \quad \mathbf{X}\in \mathcal{X} \\
& \quad \eta_{\omega} \geq 0,  \quad \forall\ \omega=1,\dots,\Omega\\
& \quad f(\mathbf{Y_\omega},\mathbf{X},\boldsymbol{\xi}_\omega) - \nu - \eta_{\omega} \leq 0,  \quad \forall\  \omega=1,\dots,\Omega  \\
& \quad \mathbf{Y}_{\omega}\in\mathcal{Y}(\mathbf{X},\boldsymbol{\xi}_{\omega}), \quad \forall\ \omega=1,\dots,\Omega   
\end{align}
\end{subequations}

However, problem (\ref{eq:SAA_two_stage_CVaR}) shares the same computational issues as the standard, expected value-based SAA formulation (\ref{eq:SAA_two_stage}), becoming intractable for a large set of scenarios. Moreover, for the general case where $f(\mathbf{Y},\mathbf{X},\boldsymbol{\xi})$ is not convex with respect to $\mathbf{X}$, there is no manageable equivalent SAA version of the problem (\ref{eq:two_stage_risk}).

\subsection{Quantile Neural Networks}
\label{sec:back_qnn}

Quantile regression introduces a probabilistic perspective on the statistical modeling of conditional variables \citep{hao2007quantile}. Let $y$ represent the response variable and
$\mathbf{X}$ the matrix of predictor variables. In the classical neural-network regression framework, the goal is to model the conditional mean $\mathbb{E}\left[ y|\mathbf{X} \right]$, typically achieved by minimizing the sum of squared residuals. The resulting model provides insights regarding the central tendency of the response variable.

On the contrary, quantile regression focuses on estimating the $\tau$-th quantile of the conditional distribution, i.e., $Q_{\tau}\left( y|\mathbf{X} \right) = \text{inf} \{y: F_{y|\mathbf{X}}(y) \geq \tau\}$ , where $F(\cdot)$ is the function to be estimated and $\tau$ lies in the interval $\left[0,1\right]$. In this case, the model is fitted by minimizing the quantile loss (sometimes called the ``pinball'' loss), which can be defined as the asymmetrically weighted sum of absolute deviations. Therefore, for a given set of $N$ data points, the quantile regression problem is set as follows:
\begin{equation}
\label{eq:quant_loss}
\mathbf{\theta} = \argmin \frac{1}{N} \sum_{i=1}^{N} \left[ \tau \epsilon_i I_{\epsilon_i \geq 0} + (1-\tau) \epsilon_i I_{\epsilon_i < 0} \right]
\end{equation}

\noindent where $\epsilon_i = y_i - f_{\theta}(\mathbf{X}_i)$ is the estimation error obtained with the model $f_{\theta}(\cdot)$ in sample $i$, $\mathbf{\theta}$ is the model weight matrix, and $I_{*}$ is an indicator function that takes the value of $1$ if the selected condition holds. As can be seen, for a quantile level $\tau$ less than $0.5$, the loss is greater when $\epsilon_i$ is negative (the prediction $f_{\theta}(\mathbf{X}_i)$ is above the actual value $y_i$). On the other hand, for a quantile level $\tau$ greater than $0.5$, the loss is greater when the error is positive. The quantile loss is only symmetrical for a $\tau$ value of $0.5$.

This probabilistic regression approach has proven to be a powerful tool for modeling and analyzing the conditional distribution of response variables, particularly in scenarios where data exhibit heteroscedasticity, non-normality, or asymmetry. Quantile neural networks (QNNs) extend this methodology into the neural network paradigm. A QNN is essentially a neural network architecture that incorporates the quantile loss function, enabling it to simultaneously estimate multiple quantiles of the conditional distribution \citep{xu2017composite}. The network is trained to minimize the weighted sum of absolute deviations for each specified quantile level, facilitating a flexible and data-driven approach to capture the variation of the response distribution.

In practice, when aiming to estimate $K$ different conditional quantiles, the problem is defined as finding network weights $\mathbf{\theta}$ that minimize the mean quantile loss across the dataset and quantile levels. This problem can be specifically expressed as: 
\begin{equation}
\label{eq:quant_loss_multi}
\mathbf{\theta} = \argmin \frac{1}{NK} \sum_{j=1}^{K} \sum_{i=1}^{N} \left[ \tau_j \epsilon_{i,j} I_{\epsilon_{i,j} \geq 0} + (1-\tau_j) \epsilon_{i,j} I_{\epsilon_{i,j} < 0} \right]
\end{equation}
where $N$ is the number of samples in the training dataset. 
In contrast to the single-quantile case (\ref{eq:quant_loss}), $\epsilon_{i,j} = y_i - f^j_{\theta}(\mathbf{X}_i)$ is now the estimation error obtained with the model on sample $i$ when predicting the $j$-th conditional quantile.

The structure of the QNN itself can be compared to the standard structure of fully connected networks in the literature. First, an input layer receives the input predictors. Next, one or more hidden layers and a final output layer can be implemented to obtain different conditional quantiles. Note that the hidden layers are arranged sequentially: a fully connected layer, followed by a nonlinear activation layer (ELU, ReLU, Sigmoid, Tanh, etc.) and possibly a dropout layer to avoid overfitting. On the other hand, the output layer is able to produce multiple quantiles that will enable the building of the dependent variable distribution.

Mathematically, an $L$-layer quantile neural network can be defined as follows. Let $\mathbf{X} \in \mathbb{R}^{n_x}$ denote the input to the network, where $n_x$ is the number of features. Each layer $l \in \{1,...,L\}$ has a weight matrix $W^{[l]}$ of dimensions $(n^{[l]}, n^{[l-1]})$, where $n^{[l]}$ is the number of neurons in layer $l$ ($n^{[0]}$ is the dimensionality of the inputs). The bias vector for layer $l$ is denoted as $b^{[l]}$ and has dimensions $(n^{[l]}, 1)$. The activation function for the hidden layers is denoted as $g^{[l]}(\cdot): \mathbb{R} \rightarrow \mathbb{R}$. The output of the $l$-th layer, denoted as $a^{[l]}$, is obtained as through the following operations:
\begin{subequations}
\begin{align}
z^{[l]} &= W^{[l]}a^{[l-1]} + b^{[l]} \\
a^{[l]} &= g^{[l]}(z^{[l]})
\end{align}
\end{subequations}

\noindent where $a^{[0]} = \mathbf{X}$ is the input to the network, and $a^{[L]}$ is the final output. The network output can be represented as $\mathbf{q} = a^{[L]}$, with $\mathbf{q} \in \mathbb{R}^{n_k}$, where $n_k$ is the number of estimated conditional quantiles.

\subsection{Mixed-Integer Formulations for ReLU Neural Networks}
\label{sec:back_embedding}

When the activation functions $g^{[l]}$ in every layer are linear or piecewise linear (e.g., ReLU, leaky ReLU), the neural network can be embedded as constraints of a mixed-integer linear program. 
Specifically, after the neural network is trained, the parameters $W^{[l]}$ and $b^{[l]}, l \in \{1,...,L\}$ are fixed, and the neural network merely denotes a learned function. 
This overall function is (piecewise) linear if all $g^{[l]}$ are (piecewise) linear, and the neural network can be reformulated by considering each activation $g(\cdot)$ separately.  
We refer the interested reader to \citet{huchette2023deep} for a comprehensive overview of this area. 

A popular method to formulate disjunctive constraints in mixed-integer programming is the so-called big-$M$ method~\citep{bonami2015mathematical}. Big-$M$ formulations are preferred owing to their simplicity and compactness, but their linear relaxations can be weak, slowing the performance of a mixed-integer solution algorithm. Therefore, stronger formulations~\citep{anderson2020strong,tsay2021partition} may be preferred for more difficult problems. Consider a single ReLU activation, i.e., a layer $l$ with $n^{[l]}=1$:
\begin{subequations}
\begin{align}
a^{[l]} &= \mathrm{max} \{ 0, W^{[l]}a^{[l-1]} + b^{[l]}  \}
\end{align}
\end{subequations}

An early set of works~\citep{fischetti2018deep,lomuscio2017approach,tjeng2018evaluating} showed that the big-$M$ method can be used to produce the following mixed-integer formulation:
\begin{subequations} \label{eqn:relu-big-m}
\begin{align}
    a^{[l]} &\geq W^{[l]}a^{[l-1]} + b^{[l]} \\
    a^{[l]} &\leq W^{[l]}a^{[l-1]} + b^{[l]} - M^-(1-\sigma)\\
    0 &\leq a^{[l]} \leq M^+ \sigma\\
    \sigma &\in \{0,1\}
\end{align}
\end{subequations}
Here, $\sigma$ is an auxiliary binary variable, while $M^+$ and $M^-$ are big-$M$ constants, which must satisfy the following bounds:
\begin{align}
    M^- \leq W^{[l]}a^{[l-1]} + b^{[l]}  \leq M^+ \label{eq:bigMbounds}
\end{align}

\newcommand{\ubar}[1]{\text{\b{$#1$}}}
Given bounds for each input variable, $a^{[l-1]}_i \in [\ubar{a}^{[l-1]}_i, \bar{a}^{[l-1]}_i ], \forall i \in \{ 1,...,n^{[l]} \}$, interval arithmetic can be used to derive valid bounds:
\begin{align}
    M^- &= \sum_i \left( \ubar{a}^{[l-1]}_i \mathrm{max}(0, W^{[l]}_i) + \bar{a}^{[l-1]} \mathrm{min}(0, W^{[l]}_i) \right) + b^{[l]} \label{eq:interval1} \\
    M^+ &= \sum_i \left( \bar{a}^{[l-1]}_i \mathrm{max}(0, W^{[l]}_i) + \ubar{a}^{[l-1]}_i \mathrm{min}(0, W^{[l]}_i) \right) + b^{[l]} \label{eq:interval2}
\end{align}

Although tighter bounds can be derived, e.g., using optimization-based bounds tightening, interval bounds are often preferred for their simplicity. 
Recent works~\citep{badilla2023computational,zhao2024bound} investigate the computational tradeoffs of more expensive bounds-tightening techniques, finding for example that interval bounds perform relatively well for shorter (less deep) neural networks. 
Several software tools, such as \texttt{JANOS}~\citep{bergman2022janos}, \texttt{MeLON}~\citep{schweidtmann2019deterministic}, and \texttt{OMLT}~\citep{ceccon2022omlt}, enable automatic translation from trained neural networks into corresponding optimization formulations, such as \eqref{eqn:relu-big-m}.

\section{Methodology}
\label{sec:methodology}

\subsection{Quantile Neural Network Methodology}
\label{sec:meth_qnn}

The key idea of our proposed methodology is to employ a QNN as a surrogate model for the second stage of two-stage stochastic optimization problems. As mentioned beforehand, the curse of dimensionality in second-stage variables and constraints incurred by the SAA methodology is the main issue for the efficient solving of two-stage stochastic problems. We aim to obtain solutions fast with a quality close to the original SAA. Furthermore, the proposed framework will be flexible: as the QNN models the distribution of the second-stage objective, it will not be limited to work with the expected value of the second stage, but can be also adapted to obtain risk-averse formulations. 

\begin{figure}[!ht]
    \centering
    \includegraphics[width=0.8\textwidth]{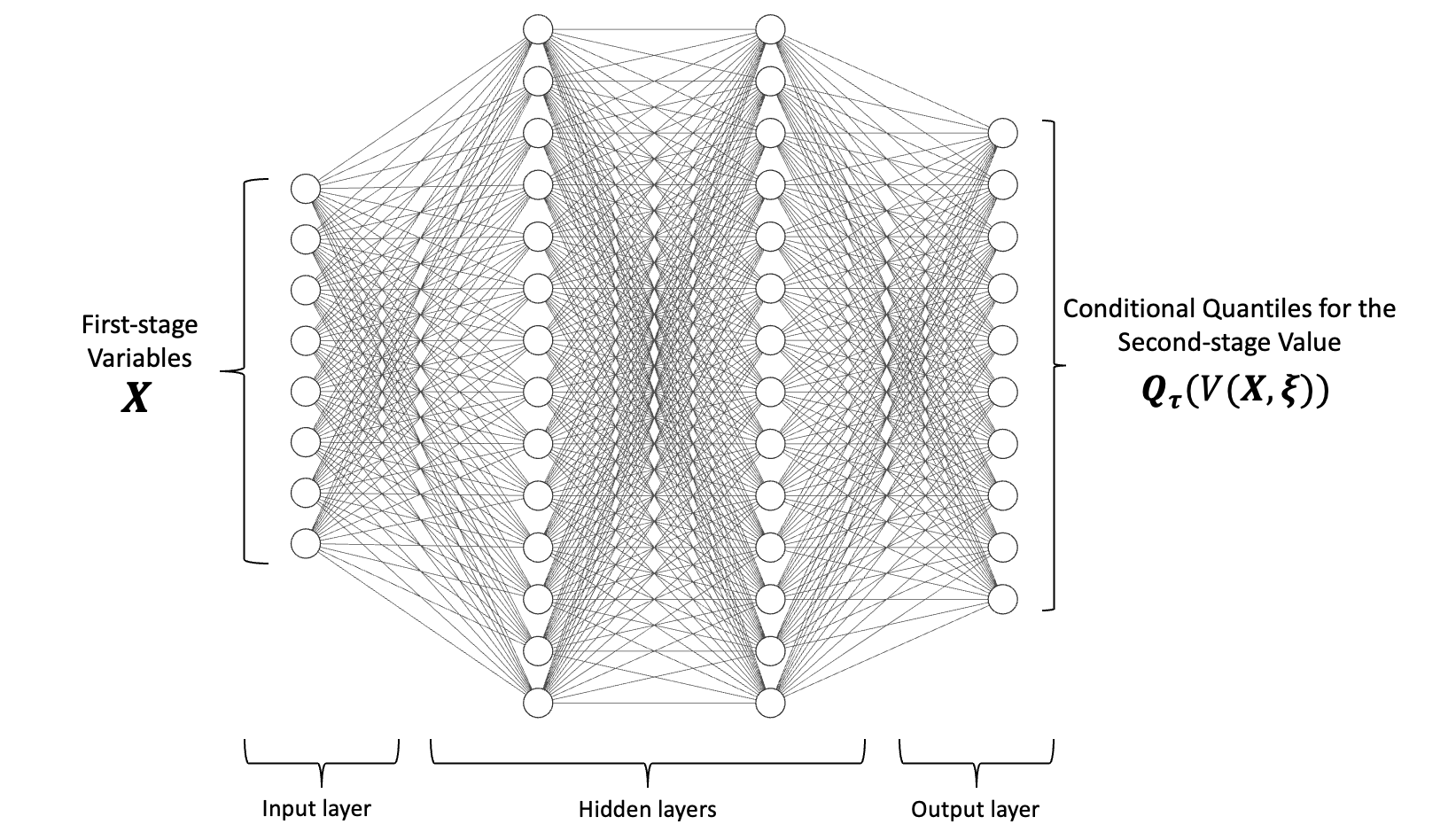}
    \caption{QNN Structure for Two-stage Optimization Problems.}
    \label{fig:qnn_2sp}
\end{figure}

With this purpose in mind, a QNN will be trained to estimate the distribution of the second-stage value function given our first-stage decision variables. In this sense, we aim to generalize the solution and the second-stage conditional distribution so that it does not depend on the number of scenarios when solving the problem. Therefore, decision variables $\mathbf{X}$ will be treated as neural network inputs, whereas the output layer will produce multiple quantiles that will enable the reconstruction of the second-stage value distribution. The structure of the QNN can be seen in Figure \ref{fig:qnn_2sp}.

In short, the QNN will learn the mapping $\mathbf{X} \rightarrow \mathbf{Q_{\tau}}(V(\mathbf{X}, \boldsymbol{\xi}))$. Unlike other methodologies that make use of neural network surrogates for two-stage stochastic programming (e.g., see \cite{patel2022neur2sp}), scenarios are not considered either as QNN inputs or in the optimization problem itself. This is achieved by learning the set of conditional quantiles $\mathbf{Q_{\tau}}(\cdot)$ instead of just estimating the mean, which can effectively characterize the second-stage value distribution $V(\mathbf{X}, \boldsymbol{\xi})$ without employing the scenarios, and allow us to account for the uncertainty directly in the output of the neural network.

Figure~\ref{fig:qnn_2sp} shows that a QNN has essentially the architecture of a multi-output, feedforward neural network (the multiple outputs are trained to correspond to quantiles of a distribution). 
As described in Section~\ref{sec:back_embedding}, such a neural network can be embedded within an optimization problem as a set of piecewise-linear constraints if all activation functions used are themselves piecewise linear.
Therefore, we use the ReLU activation function for all nodes of the QNN, such that the QNN approximation of the second-stage problem can be embedded in the first-stage problem. 
This results in a monolithic approximation of the two-stage stochastic program. 

Once the QNN surrogate model is embedded, the output layer predicts a distribution and gives the option to optimize according to different metrics. Choosing to optimize over the mean of the different outputs, that is, the conditional distribution, will create a surrogate problem for the second-stage expected value. On the other hand, deciding to optimize the conditional distribution tail mean will be equivalent to a surrogate model for risk-averse decisions. 

Some post-hoc adjustments of the QNN were considered during the development of the framework, such as conformalizing the resulting quantiles \citep{romano2019conformalized}. Conformalization methodologies for conditional quantiles aim to calibrate their coverage under finite sample guarantees. Even if this could improve the predictive performance of the adjusted QNN, it will not change the prescriptive power, as the calibration is independent of the first-stage input variables, and therefore, we do not consider it here.

However, one issue that may actually affect the solution to the problem achieved with the QNN is the ``quantile crossing'' phenomenon. Specifically, this corresponds to potential violations of the non-decreasing property of conditional quantile estimations. That is, for $\tau_1 < \tau_2 \in \left[0,1\right]$, the condition $\mathbf{Q_{\tau_1}}(V(\mathbf{X}, \boldsymbol{\xi})) < \mathbf{Q_{\tau_2}}(V(\mathbf{X}, \boldsymbol{\xi}))$ should hold for all $\mathbf{X}$. In practice, this property is challenging to guarantee during training, and recent literature has tried to address this issue in different ways. \cite{cannon2018non} adds the quantile level $\tau$ as an input feature, which will result in several QNNs if used to optimize. \cite{gasthaus2019probabilistic} propose the spline-quantile function, resulting in recurrent neural networks. In contrast, \cite{moon2021learning} develop an algorithm to achieve non-crossing QNNs, but requires using a sigmoid function on top of the hidden layers, which can not be exactly represented with a mixed-integer formulation. As can be seen, there is still lacking a straightforward approach to avoid quantile crossing for our purposes. 

\begin{figure}[!ht]
    \centering
    \includegraphics[width=0.8\textwidth]{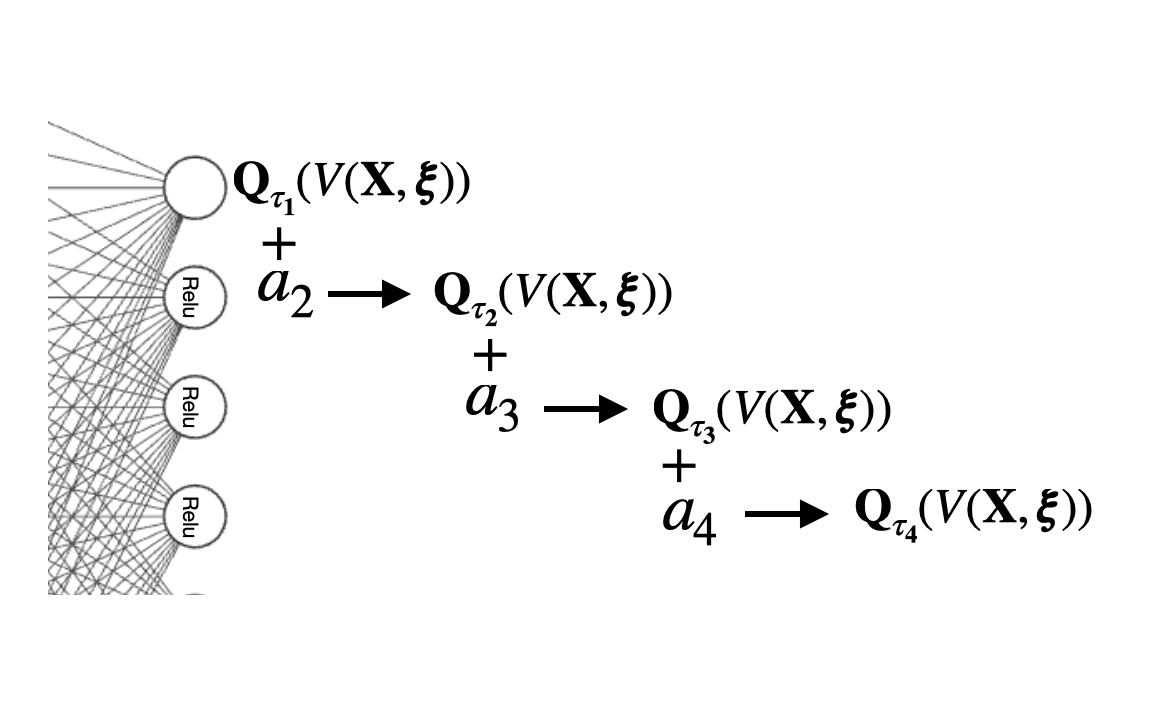}
    \caption{IQNN Output Layer Example.}
    \label{fig:iqnn_output}
\end{figure}

Nevertheless, some structural modifications can be done in the output layer to avoid the quantile-crossing phenomenon. We can predict \textit{increments} in the quantile function instead of the conditional quantiles themselves. 
In this way, quantile crossing can be avoided by using an non-negative activation function, such as ReLU. 
Then, setting the quantile estimation is easy by performing the cumulative sum until the specific $\tau$ level. This type of network will be denoted as Incremental Quantile Neural Network (IQNN), as introduced by \cite{park2022learning}. In the IQNN, decision variables $\mathbf{X}$ are still treated as input features for the network. However, ReLU activation functions are applied in the output layer (c.f. typical QNNs can employ linear activations in the output layer) to achieve an incremental behavior. This output layer is depicted in Figure \ref{fig:iqnn_output}.

In a simple way, the first neuron of the output layer will output the estimate of the lower quantile in the set, i.e., $\mathbf{Q_{\tau_1}}(V(\mathbf{X}, \boldsymbol{\xi}))$. However, from the second neuron, ReLU functions are applied, outputting a non-negative value, and allowing us to build $\mathbf{Q_{\tau_2}}(V(\mathbf{X}, \boldsymbol{\xi}))$ as the sum of $\mathbf{Q_{\tau_1}}(V(\mathbf{X}, \boldsymbol{\xi}))$ and $a_{2}$ (output of the ReLU function). This accumulation process is repeated through the nodes of the output layer.

Therefore, two neural network architectures will be considered in our framework. First, the QNN, whose structure is the classical feed-forward network with multiple outputs and ReLU activation functions in only the hidden layers. This straightforward architecture can cause quantile crossing problems, which are natural when predicting. We have tried to mimic the quantile-crossing behavior of the QNN when it has been employed as an optimization problem surrogate. 
Specifically, we can include constraints that enforce that successive quantiles must be increasing. 
Note that enforcing monotonicity in this way may lead to worse quality solutions than when we allow for some crossings, as the latter matches the setting during training. Nevertheless, allowing too much quantile crossing may lead to solutions far from the training data, leading to poor predictions. We solve this issue with additional constraints and a tolerance parameter, which will be discussed in detail in Section \ref{sec:meth_prob_formulation}. On the other hand, the IQNN solves the quantile-crossing issue by adding more non-linearities in the output layer, but allowing us not to add additional constraints or parameters.

\subsection{Data Generation}
\label{sec:meth_data}

A complete dataset is needed in order to approximate the second-stage value through the training of the (I)QNN. The data generation procedure here is fast and straightforward, as described in Algorithm \ref{alg:datagen}. As can be observed from this algorithm, only a feasible random solution $\mathbf{X}_i$ and a scenario $\boldsymbol{\xi}_i$ are needed.

\begin{algorithm}[ht]

\SetAlgoLined
\KwData{$N$: Desired size of dataset
\newline $V(\mathbf{X},\boldsymbol{\xi})$: Evaluation function for the the second-stage objective
}
\KwResult{$DF = \{(\mathbf{X}_1,v_1), ...,(\mathbf{X}_N,v_N)) \}$: Complete dataset \\ \hrulefill \\}

$DF \leftarrow \{\}$ \\
\For{i = 1 TO N}{
$\mathbf{X}_i$ $\leftarrow$ \text{Feasible random second-stage inputs} \\
$\boldsymbol{\xi}_i$ $\leftarrow$ \text{Random scenario realization} \\
Solve (\ref{second_Stage}) with fixed $\mathbf{X}_i$ and $\boldsymbol{\xi}_i$ \\
$v_i = V(\mathbf{X}_i,\boldsymbol{\xi}_i)$ \quad \text{Save second-stage value} \\
$DF \leftarrow DF \ \bigcup\ (\mathbf{X}_i, v_i)$ \\
}

\caption{Data Generation Procedure.}
\label{alg:datagen}
\end{algorithm}

With the fixed input $\mathbf{X}_i$ and the single scenario $\boldsymbol{\xi}_i$, the (second-stage only) optimization problem (\ref{second_Stage}) can be easily solved. This step may appear computationally expensive but, as the first-stage decision variables are fixed and there is a single scenario, only the recourse variables will need to be adjusted, making this problem fast and easy to solve. Finally, the value of the second-stage $v_i$ is saved, building the final dataset which is composed of ($\mathbf{X}_i, v_i$) pairs. 
Note that this data generation procedure is highly parallelizable, which enables the possibility of obtaining a decent amount of samples in just the time required to solve a second-stage problem (often on the order of seconds).

\subsection{Problem Formulation}
\label{sec:meth_prob_formulation}


We propose a surrogate adaptation of problem (\ref{eq:SAA_two_stage_CVaR}) by making use of an $L$-hidden layer ReLU (I)QNN with its corresponding mixed-integer reformulation. This surrogate model will allow us to approximate the intractable SAA, obtaining general solutions for risk-averse (or risk-neutral) two-stage stochastic optimization problems.

\subsubsection{QNN surrogate model}
\label{sec:qnn_prob_formulation}

We focus first on the QNN model. Let $a^{l}_j$ represent the ReLU output of the $j$-th neuron in the $l$-th hidden layer, for all $j \in \{1, \dots, n^l\}$ and $l \in \{1, \dots, L\}$, with $a^{0}_j$ as the $j$-th input variable of the QNN for all $j \in \{1, \dots, n^0\}$. Define $z^{l}_j$ as a continuous auxiliary variable to track the linear component (i.e., the preactivation) of the $j$-th neuron in the $l$-th layer, for all $j \in \{1, \dots, n^l\}$ and $l \in \{1, \dots, L+1\}$. Suppose $\mathbf{W}^{l}_j$ and $b^{l}_j$ denote a weight vector and a bias scalar, respectively, while $M^{-,l}_j$ and $M^{+,l}_j$ are valid big-M constants for all $j \in \{1, \dots, n^l\}$ and $l \in \{1, \dots, L\}$. The validity is defined as in \eqref{eq:bigMbounds}. 
The variable $\sigma^{l}_j \in \{0,1\}$ is binary for all $j \in \{1, \dots, n^l\}$ and $l \in \{1, \dots, L\}$, that is, for each neuron in the hidden layer. Both big-M constants and binary variables are needed for the mixed-integer reformulation of the ReLU function. Finally, let $q_{\tau_{k}}$ be the $k$-th estimated conditional quantile of the QNN at level $\tau_{k}$ for all $k \in \{1, \dots, n^k\}$ and $\Delta$ a scalar representing the tolerance towards the quantile-crossing phenomena.

With all this new set of variables, problem (\ref{eq:SAA_two_stage_CVaR}) is reformulated as surrogate problem (\ref{eq:quant2sp}), mathematically described as follows:

\begin{subequations}\label{eq:quant2sp}
\begin{alignat}{3}
\min_{\mathbf{X}} \ &(1+\lambda)c^T \mathbf{X}+ \frac{1}{n^k}\sum_{k=1}^{n^k} q_{\tau_{k}}  &+ \frac{\lambda}{(n^k - n^c)} \sum_{k=n^k - n^c}^{n^k} q_{\tau_{k}}  \qquad \qquad \qquad  \label{eq:quant2sp_OF}\\
\rm{s.t.}\notag \\
& a^{0}_j = x_j  &\forall j\in \{1, \dots, n^0\} \label{eq:quant2sp_cons1}\\
& q_{\tau_{k}} = z^{L+1}_k &\forall k\in \{1, \dots, n^k\} \label{eq:quant2sp_cons2}\\
& z^{l}_j = (\mathbf{W}^{l}_j)^T \mathbf{a}^{l-1} + b^{l}_j &\forall j\in \{1, \dots, n^l\}, l \in \{1, \dots, L+1\} \label{eq:quant2sp_cons3}\\
& a^{l}_j \geq z^{l}_j &\forall j\in \{1, \dots, n^l\}, l \in \{1, \dots, L\} \label{eq:quant2sp_cons4}\\
& a^{l}_j \leq z^{l}_j - M^{-,l}_j(1-\sigma^{l}_j) &\forall j\in \{1, \dots, n^l\}, l \in \{1, \dots, L\} \label{eq:quant2sp_cons5}\\
& 0 \leq a^{l}_j \leq M^{+,l}_j \sigma^{l}_j &\forall j\in \{1, \dots, n^l\}, l \in \{1, \dots, L\} \label{eq:quant2sp_cons6}\\
& \sigma^{l}_j \in \{0,1\} &\forall j\in \{1, \dots, n^l\}, l \in \{1, \dots, L\} \label{eq:quant2sp_cons7}\\
& q_{\tau_{k}} \leq q_{\tau_{k+1}} + \Delta &\forall k\in \{1, \dots, n^k-1\} \label{eq:quant2sp_cons8}\\
& \mathbf{X} = \{x_1, \dots, x_{n^0}\} \in \mathcal{X}  & \label{eq:quant2sp_cons9}
\end{alignat}
\end{subequations}

The objective function (\ref{eq:quant2sp_OF}) in the surrogate problem is composed of three summation terms, with the first being the weighted first-stage cost as in the original problem. The second term denotes the mean value of the conditional quantiles, which would be equivalent to the expected value of the second stage. Finally, the third term approximates the CVaR of the second stage by computing the mean of the distribution's right-hand tail, i.e., taking the mean from the $q_{\tau_{n^k-n^c}}$ conditional quantile. Observe that $\lambda$ still represents the trade-off between expectation and risk and that the minimization problem can be easily transformed into a maximization one by selecting the left tail of the distribution (i.e., worst-case scenarios for the profit).

Constraints (\ref{eq:quant2sp_cons1}) and (\ref{eq:quant2sp_cons2}) keep track of the input vector of the neural network, or the first-stage decision variables, and the output of the QNN (conditional quantiles), respectively. On the other hand, constraint (\ref{eq:quant2sp_cons1}) assigns the input of each neuron to $z_j^l$.
Constraints (\ref{eq:quant2sp_cons4}-\ref{eq:quant2sp_cons7}) represent the mixed-integer formulation of the ReLU function, as described in Section \ref{sec:back_embedding}. Constraint (\ref{eq:quant2sp_cons8}) sets the tolerance parameter $\Delta$ for the quantile crossing phenomena. Notice that each quantile $q_{\tau_{k}}$ represents a different level $\tau$, and therefore $q_{\tau_{k}}$ should be lower than $q_{\tau_{k+1}}$, as set during the training process. As this fact is not guaranteed in practice, we allow some crossing between different quantiles to mimic the QNN behavior in a prediction set-up. Finally, constraint (\ref{eq:quant2sp_cons9}) limits first-stage variables $\mathbf{X}$ to belong to the feasible set $\mathcal{X}$.

\begin{algorithm}[t]

\SetAlgoLined
\KwData{$\boldsymbol{\Delta}$: Set of tolerance parameters for evaluation
\newline $\boldsymbol{\xi}$: Stochastic scenario set
}
\KwResult{Optimal $\Delta^*$ \\ \hrulefill \\}

$\mathbf{F} \leftarrow \{\}$ \\
\For{$\Delta_i$ IN $\boldsymbol{\Delta}$}{
Solve opt. problem (\ref{eq:quant2sp}) with fixed $\Delta_i$ \\
$\mathbf{X}_i$ $\leftarrow$ Optimal solution of problem (\ref{eq:quant2sp})  \\
Solve SAA problem (\ref{eq:SAA_two_stage_CVaR}) with fixed $\mathbf{X}_i$ and $\boldsymbol{\xi}$ \\
$F_i$ $\leftarrow$ Problem (\ref{eq:SAA_two_stage_CVaR}) objective value \\
$\mathbf{F} \leftarrow \mathbf{F} \ \bigcup\ (F_i)$ \\
} 
$\Delta^* = \argmin{\mathbf{F}}$ 
\caption{Prescriptive Selection of $\Delta$ parameter}
\label{alg:tolerance_selection}
\end{algorithm}

For the optimal adjustment of parameter $\Delta$, we propose a prescriptive selection of the parameter. We describe the procedure in Algorithm \ref{alg:tolerance_selection}. For that, we built a set $\boldsymbol{\Delta}$ of values to evaluate with respect to their prescriptive performance. For each parameter in the set, the surrogate problem (\ref{eq:quant2sp}) is solved, and its solution is evaluated within an SAA problem with a fixed (small) scenario set $\boldsymbol{\xi}$. The optimal tolerance parameter $\Delta^*$ is set as the one that minimizes (or maximizes) the SAA objective function.

The developed surrogate formulation allows us to obtain solutions at different risk-aversion levels by modifying the value of $\lambda$. Setting this parameter to zero establishes the surrogate modeling of the risk-neutral problem (\ref{eq:SAA_two_stage}), while higher values of $\lambda$ will give more importance to the CVaR and thus be more risk-averse.

\subsubsection{IQNN surrogate model}
\label{sec:iqnn_prob_formulation}

Now we focus on the surrogate problem employing the IQNN model. We keep the same notation as for the QNN model (Section \ref{sec:qnn_prob_formulation}). Therefore, the resulting surrogate problem for the IQNN is defined as follows:

\begin{subequations}\label{eq:iquant2sp}
\begin{alignat}{3}
\min_{\mathbf{X}} \ &(1+\lambda)c^T \mathbf{X}+ \frac{1}{n^k}\sum_{k=1}^{n^k} q_{\tau_{k}}  &+ \frac{\lambda}{(n^k - n^c)} \sum_{k=n^k - n^c}^{n^k} q_{\tau_{k}}  \qquad \qquad \qquad  \label{eq:iquant2sp_OF}\\
\rm{s.t.}\notag \\
& a^{0}_j = x_j  &\forall j\in \{1, \dots, n^0\} \label{eq:iquant2sp_cons1}\\
& q_{\tau_{1}} = z^{L+1}_1 \label{eq:iquant2sp_cons2a}\\
& q_{\tau_{k}} = q_{\tau_{k-1}} + a^{L+1}_k &\forall k\in \{2, \dots, n^k\} \label{eq:iquant2sp_cons2b}\\
& z^{l}_j = (\mathbf{W}^{l}_j)^T \mathbf{a}^{l-1} + b^{l}_j &\forall j\in \{1, \dots, n^l\}, l \in \{1, \dots, L+1\} \label{eq:iquant2sp_cons3}\\
& a^{l}_j \geq z^{l}_j &\forall j\in \{1, \dots, n^l\}, l \in \{1, \dots, L+1\} \label{eq:iquant2sp_cons4}\\
& a^{l}_j \leq z^{l}_j - M^{-,l}_j(1-\sigma^{l}_j) &\forall j\in \{1, \dots, n^l\}, l \in \{1, \dots, L+1\} \label{eq:iquant2sp_cons5}\\
& 0 \leq a^{l}_j \leq M^{+,l}_j \sigma^{l}_j &\forall j\in \{1, \dots, n^l\}, l \in \{1, \dots, L+1\} \label{eq:iquant2sp_cons6}\\
& \sigma^{l}_j \in \{0,1\} &\forall j\in \{1, \dots, n^l\}, l \in \{1, \dots, L+1\} \label{eq:iquant2sp_cons7}\\
& \mathbf{X} = \{x_1, \dots, x_{n^0}\} \in \mathcal{X}  & \label{eq:iquant2sp_cons9}
\end{alignat}
\end{subequations}

The formulation \eqref{eq:iquant2sp} is similar to the previous one \eqref{eq:quant2sp}. The main difference is that ReLU activation functions are applied until the output layer ($L+1$), and not only in the hidden layers; see constraints (\ref{eq:iquant2sp_cons3})-(\ref{eq:iquant2sp_cons7}). Furthermore, the estimated quantiles are set differently. Only the first quantile is estimated through the linear output of the neuron in constraint (\ref{eq:iquant2sp_cons2a}), that is, the ReLU is not applied so the output can be negative, if needed. The following quantiles are built in an incremental sense adding to the previous quantile the non-negative output of its respective ReLU function, as in constraint (\ref{eq:iquant2sp_cons2b}).
The objective function formulation and the rest of the possible constraints remain as in the aforementioned model \eqref{eq:quant2sp}, and the quantile crossing constraint is not needed by the construction of the quantile outputs, where monotonicity is ensured.

\subsection{Related Work}
\label{sec:meth_related_work}

As has been stated throughout the exposition of this paper, recent works~\citep{patel2022neur2sp,kronqvist2023alternating} have used embedded NNs as surrogate models in a general framework for two-stage stochastic optimization. 
Here we compare our approach to the presented \texttt{Neur2SP} framework~\citep{patel2022neur2sp}. \citet{kronqvist2023alternating} propose adaptive-sampling procedures for this framework. 
Their two presented NN architectures (NN-E and NN-P) focus on approximating the expected value of the second stage using first-stage decisions and scenarios as input features for the model. Table~\ref{tab:qnn_neur2sp_generalcomp} gives the main differences between our (I)QNN-based approach and these two architectures. 

\begin{table}[ht]
    \centering
    \begin{adjustbox}{width=1\textwidth}
    \begin{tabular}{lc|cc}
        \toprule
                 & \textbf{(I)QNN} & \textbf{NN-E} & \textbf{NN-P} \\
        \midrule
        Data generation & $v_i = V(\mathbf{X}_i,\boldsymbol{\xi}_i)$ & $v_i^* = \sum_s p_s V(\mathbf{X}_i,\boldsymbol{\xi}_s) $ & $v_i = V(\mathbf{X}_i,\boldsymbol{\xi}_i)$ \\
        \midrule
        Training Data & $\{(\mathbf{X}_n,v_n)\}_{n=1}^N$ & $\{(\mathbf{X}_n,\boldsymbol{\xi}_{\lambda},v_n^*)\}_{n=1}^N$ & $\{(\mathbf{X}_n,\boldsymbol{\xi}_n,v_n)\}_{n=1}^N$ \\
        \midrule
        NN training loss & Quantile loss & MSE & MSE \\
        \midrule
        Opt. cost & \makecell{Single\\NN embedding} & \makecell{Single\\NN embedding} & \makecell{Multiple\\NN embedding}  \\
        \midrule
        Allows Exp. value opt. & Yes & Yes & Yes \\
        \midrule
        Allows Risk opt. (e.g. CVaR) & Yes & No & No \\
        \bottomrule
    \end{tabular}
    \end{adjustbox}
    \caption{Complete framework comparison with the original \texttt{Neur2SP} framework~\citep{patel2022neur2sp}.}
    \label{tab:qnn_neur2sp_generalcomp}
\end{table}

Our data generation procedure is the same as the one proposed in the NN-P architecture. This procedure, with a fixed random solution and single scenario, is extremely fast and parallelizable. In contrast, the data-generation step for the NN-E architecture is its main bottleneck, as it needs to solve the problem for a complete scenario set, rather than for a unique scenario value.

NN-E and NN-P models are trained with scenarios (or embedding of scenarios) in the input layer, and their main goal is to estimate the expected value of the second stage. Therefore, both architectures aim to minimize the mean square error. Besides, one important limitation of the NN-P is that it requires one embedding (i.e., one additional NN) per scenario.

In contrast, we do not save the values of the scenarios, as they will not take part in the (I)QNN architecture. Rather, we minimize the quantile loss in order to have a proper generalization of the second-stage distribution using a single NN architecture, as is the case with NN-E. This means that our approach will not be an expert system for optimizing with the expected value, but instead will be flexible enough to obtain solutions both in the expected value and in the Conditional Value-at-Risk (CVaR) of the second stage.

\section{Case Studies}
\label{sec:cases_estudies}

In this section, we will focus on evaluating the performance of our (I)QNN-based framework within a wide range of two-stage stochastic optimization problems and benchmarks. Firstly, we introduce the optimization problems and the networks employed as surrogate models. Then, we evaluate the impact of the dataset size in terms of the (I)QNN predictive and prescriptive degradation. Furthermore, we compare optimization solution times and the true objective obtained by the quantile methodology with the ones obtained by NN-E and NN-P (in risk-neutral optimization), and SAA (in both risk-neutral and risk-averse optimization). Finally, we analyze the quantile crossing phenomena and its prescriptive impact in terms of solution quality.

\subsection{Experimental setup}
For all experiments, a desktop workstation with Intel Core i7 11700 CPU, 64 GB RAM, and a NVIDIA GeForce GTX 2060 graphics card was employed. Results were obtained using Pytorch 2.1 \citep{paszke2019pytorch} for model training and Gurobi 11.0 \citep{gurobi} as the optimization solver.

\subsubsection{Two-stage stochastic optimization problems}

For the sake of comparative analysis, we adopt a collection of two-stage stochastic optimization problems previously explored by \cite{patel2022neur2sp}. This set includes:
\begin{itemize}
\item Three variants of a Capacitated Facility Location Problem  \citep{cornuejols1991comparison}, a linearly constrained minimization problem with binary variables in both the first and second stages. These problems are denoted as CFLP-n-m, where n is the number of facilities and m is the number of customers.
\item A maximization Investment Problem \citep{schultz1998solving} (IP-I-H) that includes continuous variables in the first stage and binaries in the second, with linear constraints in both stages. 
\end{itemize}

A detailed description of these problems can be found in \cite{patel2022neur2sp}.
In total, we consider four two-stage stochastic optimization problems, and evaluate them across varying numbers of scenarios. Furthermore, while originally formulated for a risk-neutral approach, we have introduced modifications to the optimization problems to enable the exploration of risk-averse solutions within their SAA form, so we can also evaluate the performance of our (I)QNN framework in risk-averse settings.

\subsubsection{(I)QNN model selection}

In the pursuit of optimal hyperparameter values for training (I)QNNs, we perform a grid-search across a pool of 100 potential configurations, establishing different models for each optimization problem. This process is highly parallelizable and relatively fast, as compact networks can give good results when used as surrogate models. For the rest of this section, we set the output layer of every network to have 50 neurons, which represents 50 equally-spaced quantiles representing $\tau$ from 0.01 to 0.99. We find this number of quantiles to generally provide a good representation of the distribution function.
The selection of the most favorable configuration is based on minimizing the validation quantile loss, with details presented in \ref{appendix:model_sel}. The validation set is constructed through a 20\% partition of the complete dataset. Throughout the training phases of all networks, a dataset comprising 20,000 samples is generated. We study the importance of this sample size selection in Section \ref{sec:data_size}.

\subsection{Dataset size impact}
\label{sec:data_size}

Before delving into the performance analysis of the developed framework, we study the impact that the data generation procedure has on the results obtained by the surrogate models. With that purpose in mind, we train both QNN and IQNN models with different numbers of training samples and evaluate their validation quantile losses and the true objectives that we would obtain using the heuristic decisions by evaluating them in several scenario sets.

\begin{figure}[ht]
    \centering
    \begin{subfigure}{}
    \centering
    \includegraphics[width=0.48\textwidth]{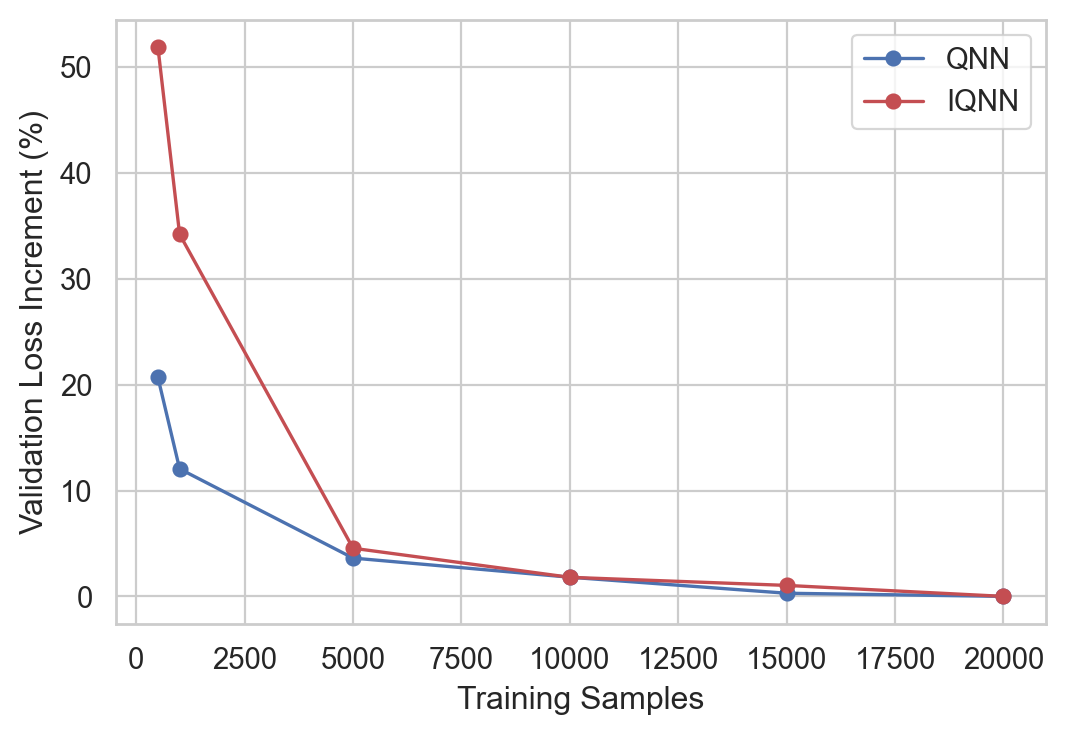}
    \end{subfigure}
    \begin{subfigure}{}
    \centering
    \includegraphics[width=0.48\textwidth]{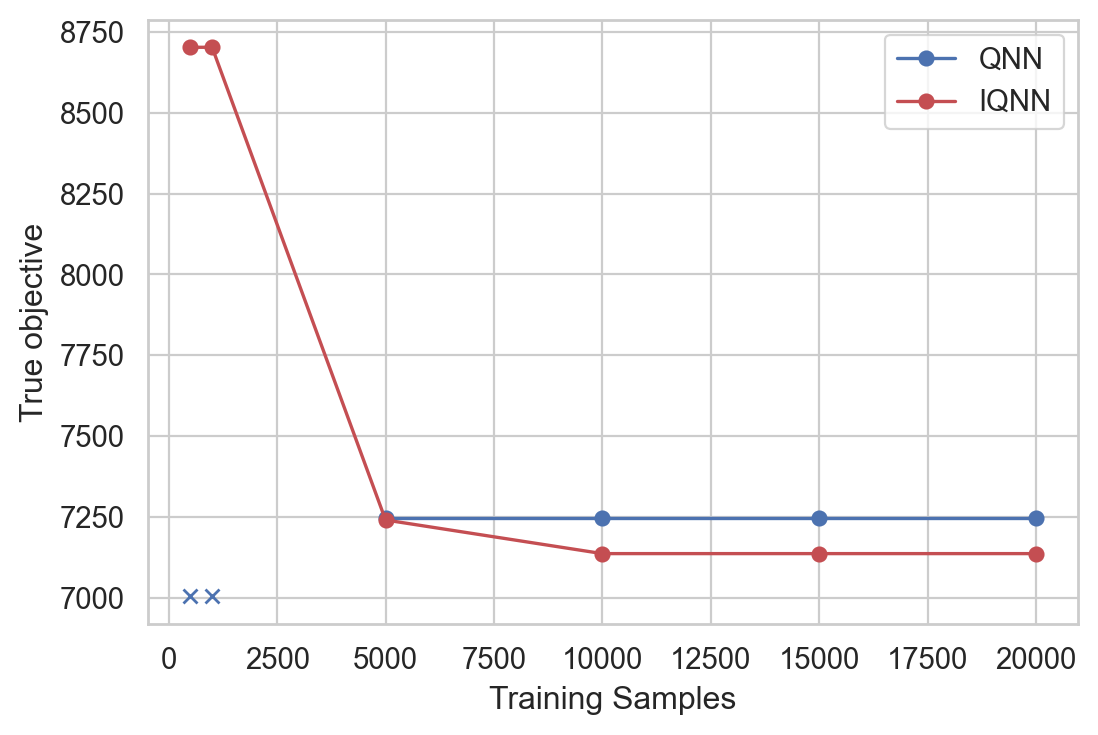}
    \end{subfigure}
    \caption{Predictive (left) and prescriptive (right) performance of the (I)QNN framework with different numbers of training samples for the risk-averse CFLP-10-10 problem.}
    \label{fig:loss_obj_size}
\end{figure}

CFLP-10-10 is chosen as an illustrative optimization problem for this section with a risk-neutral objective. QNN and IQNN models are then trained with different numbers of samples in the training phase, which is done after a hyper-parameter tuning procedure. 
Figure \ref{fig:loss_obj_size} shows the performance of the QNN and IQNN models with varying number of training samples. On the left of the figure, we see the validation loss increase (\%) compared to the model trained with all 20,000 samples. For both QNN and IQNN models, diminishing returns in model validation accuracy can be observed, with the validation loss increment close to zero for $\geq$10,000 samples.

We can observe the optimal objective (cost) of the corresponding heuristic solutions in Figure \ref{fig:loss_obj_size} (right). Specifically, we show the mean value of the true objective evaluated for 10 random scenario sets of size 500. In general, the objective value found appears stable when using $\geq$5,000 samples for training, for both the both QNN and IQNN. Note that, for the QNN model, the best true objective appears to be obtained by training the model with fewer than 1,000 samples. However, this is likely an artifact of the approximation error---we face a high number of infeasibilities that are only avoided by setting the quantile-crossing parameter $\Delta$ to a large value of 500.

To conclude, we observe that the framework performs well when training the networks with at least 10,000 samples, which can be obtained quickly with the parallelizable data generation procedure (Algorithm \ref{alg:datagen}). However, we set the number of training samples in our experiments to 20,000, in anticipation of more complex relationships to learn.

\subsection{Results for risk-neutral optimization}

We study the (I)QNN methodology's ability to obtain high-quality heuristic risk-neutral solutions, that is, to optimize the expected value of the second stage. For that, we compare its performance with the existing \texttt{Neur2SP} framework~\citep{patel2022neur2sp}, which acts as an expert system for risk-neutral optimization. In addition, the standard SAA is also considered, as it is the general approach that can be used for the vast majority of problem structures.
Note that all of the above approaches serve as approximations to the general two-stage stochastic program, as discussed in Section~\ref{sec:background}.

First, we report the total computational times to formulate (including data generation and model training) and solve the benchmark problems using the proposed quantile-based framework, compared with the rest of the competitive approaches in Table \ref{tab:times_risk_neutral}. For comparison purposes, we reproduce the total times for NN-E, NN-P, and SAA  reported in \cite{patel2022neur2sp}, where the time limit for SAA is set to three hours. The total times reported for NN-E and NN-P include data generation, training of the network, and optimization of the combined problem. The times for the latter optimization step specifically are shown in parentheses.
Note that the times for NN-P and SAA are dependent on the number of training samples used. 

\begin{table}
    \centering
    \begin{adjustbox}{width=1\textwidth}
    \begin{tabular}{lrrrrr|rrrr|rrr}
        \toprule
         & \multicolumn{5}{c}{\textbf{QNN}} & \multicolumn{4}{c}{\textbf{IQNN}} & \multicolumn{1}{c}{\textbf{NN-E}} & \multicolumn{1}{c}{\textbf{NN-P}} & \multicolumn{1}{c}{\textbf{SAA}} \\
        \cmidrule(lr){2-13}
        Problem & \thead{Data\\Generation} & Training & \thead{Tolerance\\Selection} & \thead{Problem\\Solving} & \thead{\textbf{Total}\\\textbf{Time}} & \thead{Data\\Generation} & Training & \thead{Problem\\Solving} & \thead{\textbf{Total}\\\textbf{Time}} &\thead{\textbf{Total}\\\textbf{Time}} & \thead{\textbf{Total}\\\textbf{Time}} & \thead{\textbf{Total}\\\textbf{Time}} \\
        \midrule
        CFLP-10-10-100sc & \multirow{3}{*}{38.83} & \multirow{3}{*}{526.43} & \multirow{3}{*}{16.40} &\multirow{3}{*}{3.30} & \multirow{3}{*}{584.96} & \multirow{3}{*}{38.83} & \multirow{3}{*}{407.01} & \multirow{3}{*}{0.04} & \multirow{3}{*}{445.88} & 2,490.73 (0.38) & 148.99 (8.28) & 4,410.60 \\
        CFLP-10-10-500sc &                        &                         &  &                     &                       & & & & & 2,490.95 (0.60) & 347.01 (206.30) & 10,800.17 \\
        CFLP-10-10-1000sc &                       &                         &   &                    &                       & & &  &  & 2,490.99 (0.64) & 997.47 (856.77) & 10,800.87 \\
        \midrule
        CFLP-25-25-100sc & \multirow{3}{*}{385.14} & \multirow{3}{*}{417.60} & \multirow{3}{*}{47.30} &\multirow{3}{*}{0.32} & \multirow{3}{*}{850.36} & \multirow{3}{*}{385.14} & \multirow{3}{*}{383.92} & \multirow{3}{*}{0.03} & \multirow{3}{*}{769.09}  & 6,354.50 (0.44) & 957.76 (4.86) & 10,800.06 \\
        CFLP-25-25-500sc &                        &                         &   &                    &                      & & & &   &  6,354.60 (0.54) & 979.31 (26.41) & 10,800.14 \\
        CFLP-25-25-1000sc &                       &                         &   &                    &                       & & & &  &  6,354.64 (0.58) & 1,007.35 (54.45) & 10,800.36 \\
        \midrule
        CFLP-50-50-100sc & \multirow{3}{*}{725.31} & \multirow{3}{*}{283.43} & \multirow{3}{*}{207.30} & \multirow{3}{*}{0.32} & \multirow{3}{*}{1,216.36} & \multirow{3}{*}{725.31} & \multirow{3}{*}{258.17} & \multirow{3}{*}{0.03} & \multirow{3}{*}{1,010.51} & 8,163.28 (1.66) & 284.78 (21.10) & 10,800.05 \\
        CFLP-50-50-500sc &                        &                         &   &                    &                       &&&&  &  8,162.87 (1.25) & 437.31 (173.63) & 10,806.15 \\
        CFLP-50-50-1000sc &                       &                         &    &                   &                       &&&&  &  8,163.06 (1.44) & 835.80 (572.12) & 10,805.82 \\
        \midrule
        IP-I-H-441sc & \multirow{3}{*}{15.46} & \multirow{3}{*}{259.28} & \multirow{3}{*}{4.20} & \multirow{3}{*}{0.06} & \multirow{3}{*}{279.00} & \multirow{3}{*}{15.46} & \multirow{3}{*}{270.00} & \multirow{3}{*}{0.09} & \multirow{3}{*}{285.55} & 9,338.79 (0.32) & 1,409.06 (1,231.48) & 10,800.00 \\
        IP-I-H-1681sc &                        &                         &  &                     &                       &&&&  &  9,338.80 (0.33) & 10,994.47 (10,816.89) & 10,800.03 \\
        IP-I-H-10000sc &                       &                         &  &                     &                        &&&& &  9,338.85 (0.38) &    ---    & 10,802.10 \\
        \bottomrule
    \end{tabular}
    \end{adjustbox}
    \caption{Computing times (in seconds) for the (I)QNN framework and competitive approaches in risk-neutral optimization. NN-E, NN-P, and SAA times are reproduced from \cite{patel2022neur2sp}. CPU times for solving the resulting optimization problems for NN-E and NN-P are shown in parentheses.}
    \label{tab:times_risk_neutral}
\end{table}

In particular, it is important to note that the proposed quantile-based approaches produce general solutions to the optimization problems, as the scenario set does not take part in the resulting surrogate model-based formulation. Similarly, the NN-E and NN-P approaches of \texttt{Neur2SP} can obtain solutions with different scenario-set sizes but, as in the quantile case, the network surrogate model has to be trained only once.

In general, we observe that the quantile frameworks are extremely fast in solving the resulting optimization problem, at a similar level as NN-E (times in parentheses). This is because in both cases, a single NN is embedded in the surrogate problem, in contrast to NN-P, where one embedding per scenario is needed. Furthermore, the set-up for (I)QNN, that is, data generation and training of the network is as fast as in NN-P. Therefore, we combine the strengths of both NN-E (fast problem solution) and NN-P (fast data generation).

The problem-solving times with the QNN surrogate model range between 0.06 and 3.30 seconds, while the total methodological time is between 279 and $\sim$1,200 seconds. The problem-solving times with the IQNN surrogate model are below 0.1 seconds, with a total methodological time between 300 and 1,000 seconds. These times are generally lower than NN-E and NN-P ones, and significantly lower than the those of the SAA, where the majority of problems remain unsolved after the three-hour time limit.

To study the quality of the obtained solutions, we empirically evaluate them within different sets of scenarios. Table \ref{tab:true_obj_risk_neutral} shows the mean true objective obtained by evaluating the solution given by the surrogate model within 10 different scenario sets of a given size. Results from NN-E, NN-P, and SAA are again reproduced from \cite{patel2022neur2sp}.

\begin{table}
    \centering
    \begin{adjustbox}{width=1\textwidth}
    \begin{tabular}{lrr|rrr}
        \toprule
        Problem & \textbf{QNN} & \textbf{IQNN} & \textbf{NN-E} & \textbf{NN-P} & \textbf{SAA} \\
        \midrule
        CFLP-10-10-100sc & 7,129.68 (1.93\%) & 7,124.11 (1.85\%) & 7,174.57 & 7,109.62 & \textbf{6,994.77} \\
        CFLP-10-10-500sc &   7,136.43 (1.90\%) & 7,114.86 (1.59\%) & 7,171.79  & 7,068.91  & \textbf{7,003.30}  \\
        CFLP-10-10-1000sc &  7,108.05 (0.27\%) & 7,095.69 (0.10\%) & 7,154.60  & \textbf{7,040.70} &  7,088.56  \\
        \midrule
        CFLP-25-25-100sc & 11,967.66 (0.87\%) & 11,999.41 (1.13\%) & \textbf{11,773.01} & \textbf{11,773.01} & 11,864.83  \\
        CFLP-25-25-500sc &  11,882.88 (-2.36\%) & 11,915.06 (-2.10\%) & \textbf{11,726.34}  & \textbf{11,726.34}  & 12,170.67  \\
        CFLP-25-25-1000sc & 11,870.60 (0.02\%) & 11,915.35 (0.40\%) &  \textbf{11,709.90}  &  \textbf{11,709.90}  & 11,868.04  \\
        \midrule
        CFLP-50-50-100sc & 25,944.94 (2.35\%) & 27,603.63 (8.89\%) & 25,236.33 & \textbf{25,019.64} & 25,349.21  \\
        CFLP-50-50-500sc & 25,906.90 (-7.60\%) & 27,324.21 (-2.54\%) & 25,281.13  & \textbf{24,964.33}  & 28,037.66  \\
        CFLP-50-50-1000sc &  25,881.86 (-14.53\%) & 27,178.49 (-10.25\%) & 25,247.77  & \textbf{24,981.70}  &  30,282.41  \\
        \midrule
        IP-I-H-441sc & 65.91 (-1.98\%) & 66.36 (-1.31\%) & 65.12 & 65.12 & \textbf{67.24} \\
        IP-I-H-1681sc & 65.60 (0.29\%) & \textbf{65.74} (0.50\%) & 65.63 & 65.34 & 65.41 \\
        IP-I-H-10000sc & \textbf{65.89} (1.95\%) & 65.84 (1.87\%) & 65.66  &  ---   &  64.63 \\
        \bottomrule
    \end{tabular}
    \end{adjustbox}
    \caption{True objective results for risk-neutral optimization. Relative gaps between QNN and IQNN objective values and the SAA approach are shown in brackets. The best results are highlighted in bold. Results from NN-E, NN-P, and SAA are reproduced from~\cite{patel2022neur2sp}.}
    \label{tab:true_obj_risk_neutral}
\end{table}

Comparing the proposed quantile-based approaches with SAA, we find that, for the QNN structure, improvements of the true objective of up to 14.5\% are found, whereas in the cases where SAA outperforms QNN, the gap is no more than 2.35\%. On the other hand, relative gaps between IQNN and SAA range from a deterioration of 8.89\% to an improvement of 10.25\%. In general, the biggest improvements are seen in optimization problems with many first-stage decision variables or with a high number of scenarios considered. 

On the other hand, the \texttt{Neur2SP} framework generally produces slightly better solutions than the quantile-based framework. This is an expected result, as the framework has been developed specifically to produce risk-neutral heuristic solutions in all steps, from the data generation procedure to the training and embedding of the models. 
Nevertheless, considering both the advantages of the proposed framework in terms of computational times (Table \ref{tab:times_risk_neutral}) and quality approximations of the true objective results (Table \ref{tab:true_obj_risk_neutral}), together with its potential to model risk aversion, we observe that the proposed quantile-based framework provides a promising trade-off between computational time and the quality of its solutions.

\subsection{Results for risk-averse optimization}

We next evaluate the ability of the (I)QNN framework to produce risk-averse solutions. With that purpose in mind, we adapt the previously used two-stage optimization problems to their mean-risk formulations (\ref{eq:SAA_two_stage_CVaR}) so we can obtain the SAA solution.
In this section, the proposed quantile-based framework is only compared with the SAA. To our knowledge, no other general machine learning-based framework has been developed to tackle two-stage risk-averse optimization through model embedding. 

Table \ref{tab:time_risk_averse} shows the computing times for the optimization step of the (I)QNN framework and the SAA approach. In this case, we only show the tolerance parameter selection and the problem-solving times, as the times for data generation and model training steps are consistent with those from the previous case study (see Table \ref{tab:times_risk_neutral}). We emphasize that, once the quantile-based model is trained, the same model can be used for \textit{both} risk-neutral and risk-averse optimization formulations, including different values of the CVaR importance parameter $\lambda$ and the specific quantile under consideration $(1-\alpha)$. 
Specifically, predicting the distribution of the second-stage objective allows us to formulate generic risk-averse optimization problems. 
We set the time limit for the SAA approach to two hours.

\begin{table}
    \centering
    \begin{adjustbox}{totalheight=0.95\textheight} 
    \begin{tabular}{lccc|cc|c|r}
        \toprule
                &           &          &           & \multicolumn{2}{c}{\textbf{QNN}} & \multicolumn{1}{c}{\textbf{IQNN}}  & \multicolumn{1}{c}{\textbf{SAA}} \\
        \cmidrule(lr){2-8}
        Problem & $\lambda$ & $\alpha$ & Scenarios &  \thead{Tolerance\\Selection} & \thead{Problem\\Solving} & \thead{Problem\\Solving} & \thead{Problem\\Solving} \\
        \midrule
        \multirow{12}{*}{CFLP-10-10 (-)} & \multirow{4}{*}{0.1} & \multirow{2}{*}{0.7}  & 500 & \multirow{2}{*}{20.41} & \multirow{2}{*}{1.95} & \multirow{2}{*}{0.04} & 7,203.38\\
                                     &                      &                          & 1000 &                         &                       &                       & 7,227.03 \\
                                     \cmidrule(lr){3-8}
                                     &                      & \multirow{2}{*}{0.9}  & 500 & \multirow{2}{*}{19.66} & \multirow{2}{*}{2.17} & \multirow{2}{*}{0.04} & 7,205.36 \\
                                     &                      &                       & 1000 &                         &                       &                     & 7,200.57 \\
                                     \cmidrule(lr){2-8}
                                     & \multirow{4}{*}{0.5} & \multirow{2}{*}{0.7}  & 500  & \multirow{2}{*}{21.03} & \multirow{2}{*}{2.02} & \multirow{2}{*}{0.02} & 7,200.58 \\
                                     &                      &                       & 1000 &                         &                       &                       & 7,233.79 \\
                                     \cmidrule(lr){3-8}
                                     &                      & \multirow{2}{*}{0.9}  & 500 & \multirow{2}{*}{20.37} & \multirow{2}{*}{1.98} & \multirow{2}{*}{0.04} & 1,560.71 \\
                                     &                      &                       & 1000 &                         &                       &                     & 7,227.18 \\
                                     \cmidrule(lr){2-8}
                                     & \multirow{4}{*}{1} & \multirow{2}{*}{0.7}  & 500 & \multirow{2}{*}{19.86} & \multirow{2}{*}{1.99} & \multirow{2}{*}{0.03} & 7,211.25 \\
                                     &                      &                     & 1000 &                         &                       &                      & 7,219.32 \\
                                     \cmidrule(lr){3-8}
                                     &                      & \multirow{2}{*}{0.9}  & 500 & \multirow{2}{*}{21.73} & \multirow{2}{*}{2.37} & \multirow{2}{*}{0.04} & 1,015.16 \\
                                     &                      &                     & 1000  &                         &                       &                      & 7,268.38 \\
        \midrule
        \multirow{12}{*}{CFLP-25-25 (-)} & \multirow{4}{*}{0.1} & \multirow{2}{*}{0.7}  & 500 & \multirow{2}{*}{110.30} & \multirow{2}{*}{0.20} & \multirow{2}{*}{0.02} & 7,200.39\\
                                     &                      &                          & 1000 &                         &                       &                       & 7,200.46 \\
                                     \cmidrule(lr){3-8}
                                     &                      & \multirow{2}{*}{0.9}  & 500 & \multirow{2}{*}{61.15} & \multirow{2}{*}{0.24} & \multirow{2}{*}{0.02} & 7,200.18 \\
                                     &                      &                       & 1000 &                         &                       &                     & 7,200.33 \\
                                     \cmidrule(lr){2-8}
                                     & \multirow{4}{*}{0.5} & \multirow{2}{*}{0.7}  & 500  & \multirow{2}{*}{123.45} & \multirow{2}{*}{0.23} & \multirow{2}{*}{0.03} & 7,200.46 \\
                                     &                      &                       & 1000 &                         &                       &                       & 7,201.58 \\
                                     \cmidrule(lr){3-8}
                                     &                      & \multirow{2}{*}{0.9}  & 500 & \multirow{2}{*}{70.05} & \multirow{2}{*}{0.30} & \multirow{2}{*}{0.03} & 7,200.64 \\
                                     &                      &                       & 1000 &                         &                       &                     & 7,213.66 \\
                                     \cmidrule(lr){2-8}
                                     & \multirow{4}{*}{1} & \multirow{2}{*}{0.7}  & 500 & \multirow{2}{*}{125.35} & \multirow{2}{*}{0.27} & \multirow{2}{*}{0.02} & 7,200.31 \\
                                     &                      &                     & 1000 &                         &                       &                      & 7,204.80 \\
                                     \cmidrule(lr){3-8}
                                     &                      & \multirow{2}{*}{0.9}  & 500 & \multirow{2}{*}{116.11} & \multirow{2}{*}{0.31} & \multirow{2}{*}{0.01} & 7,200.49 \\
                                     &                      &                     & 1000  &                         &                       &                      & 7,200.71 \\
        \midrule
        \multirow{12}{*}{CFLP-50-50 (-)} & \multirow{4}{*}{0.1} & \multirow{2}{*}{0.7}  & 500 & \multirow{2}{*}{169.60} & \multirow{2}{*}{5.26} & \multirow{2}{*}{0.11} & 7,224.05\\
                                     &                      &                          & 1000 &                         &                       &                       & 7,200.47 \\
                                     \cmidrule(lr){3-8}
                                     &                      & \multirow{2}{*}{0.9}  & 500 & \multirow{2}{*}{179.31} & \multirow{2}{*}{4.59} & \multirow{2}{*}{0.12} & 7,200.45 \\
                                     &                      &                       & 1000 &                         &                       &                     & 7,200.37 \\
                                     \cmidrule(lr){2-8}
                                     & \multirow{4}{*}{0.5} & \multirow{2}{*}{0.7}  & 500  & \multirow{2}{*}{200.40} & \multirow{2}{*}{3.12} & \multirow{2}{*}{0.06} & 7,202.31 \\
                                     &                      &                       & 1000 &                         &                       &                       & 7,200.41 \\
                                     \cmidrule(lr){3-8}
                                     &                      & \multirow{2}{*}{0.9}  & 500 & \multirow{2}{*}{185.77} & \multirow{2}{*}{3.25} & \multirow{2}{*}{0.04} & 7,201.73 \\
                                     &                      &                       & 1000 &                         &                       &                     & 7,200.51 \\
                                     \cmidrule(lr){2-8}
                                     & \multirow{4}{*}{1} & \multirow{2}{*}{0.7}  & 500 & \multirow{2}{*}{186.01} & \multirow{2}{*}{3.83} & \multirow{2}{*}{0.07} & 7,204.80 \\
                                     &                      &                     & 1000 &                         &                       &                      & 7,200.54 \\
                                     \cmidrule(lr){3-8}
                                     &                      & \multirow{2}{*}{0.9}  & 500 & \multirow{2}{*}{187.62} & \multirow{2}{*}{4.63} & \multirow{2}{*}{0.05} & 7,200.69 \\
                                     &                      &                     & 1000  &                         &                       &                      & 7,221.34 \\
        \midrule
        \multirow{12}{*}{IP-I-H (+)}  & \multirow{4}{*}{0.1} & \multirow{2}{*}{0.7}  & 500 & \multirow{2}{*}{4.20} & \multirow{2}{*}{0.01} & \multirow{2}{*}{0.11} & 7,200.08\\
                                     &                      &                          & 1000 &                         &                       &                       & 7,201.43 \\
                                     \cmidrule(lr){3-8}
                                     &                      & \multirow{2}{*}{0.9}  & 500 & \multirow{2}{*}{4.40} & \multirow{2}{*}{0.01} & \multirow{2}{*}{0.10} & 7,200.22 \\
                                     &                      &                       & 1000 &                         &                       &                     & 7,209.53 \\
                                     \cmidrule(lr){2-8}
                                     & \multirow{4}{*}{0.5} & \multirow{2}{*}{0.7}  & 500  & \multirow{2}{*}{4.40} & \multirow{2}{*}{0.01} & \multirow{2}{*}{0.10} & 7,201.58 \\
                                     &                      &                       & 1000 &                         &                       &                       & 7,206.25 \\
                                     \cmidrule(lr){3-8}
                                     &                      & \multirow{2}{*}{0.9}  & 500 & \multirow{2}{*}{4.40} & \multirow{2}{*}{0.01} & \multirow{2}{*}{0.10} & 7,200.19 \\
                                     &                      &                       & 1000 &                         &                       &                     & 7,243.12 \\
                                     \cmidrule(lr){2-8}
                                     & \multirow{4}{*}{1} & \multirow{2}{*}{0.7}  & 500 & \multirow{2}{*}{4.50} & \multirow{2}{*}{0.01} & \multirow{2}{*}{0.10} & 7,200.15 \\
                                     &                      &                     & 1000 &                         &                       &                      & 7,208.24 \\
                                     \cmidrule(lr){3-8}
                                     &                      & \multirow{2}{*}{0.9}  & 500 & \multirow{2}{*}{4.40} & \multirow{2}{*}{0.01} & \multirow{2}{*}{0.10} & 7,200.08 \\
                                     &                      &                     & 1000  &                         &                       &                      & 7,233.84 \\
        \bottomrule
    \end{tabular}
    \end{adjustbox}
    \caption{Computing times (in seconds) for the (I)QNN framework and SAA in risk-averse optimization.}
    \label{tab:time_risk_averse}
\end{table}

We observe that solving times for the quantile-based models are considerably lower than those for the SAA. For the QNN model, solving times range between 0.01 to slightly more than 5 seconds. The prescriptive tolerance parameter selection does not take more than 200 seconds, even for the more extensive problems. On the other hand, the IQNN-based framework requires no more than 0.12 seconds to produce a solution for the surrogate problem. In contrast, the SAA reaches the time limit in almost all cases.

To provide a better context for the real improvement made in computing times, assume we have already trained an IQNN model. Consider that, once the model is trained, we could sequentially obtain all the heuristic solutions for an optimization problem (using the same surrogate model) with three possible values of $\lambda$ and two values of $\alpha$ in less than 1 second, while the same procedure would require around 12 hours using the SAA approach.

As in the previous case study, we evaluate the quality of the heuristic solutions employing different scenario sets. Table \ref{tab:true_obj_averse} compares the true objective values when employing these heuristic solutions. Values presented for (I)QNN are the means across a different number of random evaluations (ten evaluations in the CFLP-10-10 and IP-I-H and five for the rest of the problems), while we present the best objective value of the SAA found by the end of its optimization time. Results are broken down by risk-aversion level $\lambda$, the tail of the distribution $\alpha$, and the size of the scenario set where we evaluate the solution.

\begin{table}
    \centering
    \begin{adjustbox}{width=0.8\textwidth, totalheight=0.95\textheight} 
    \begin{tabular}{lccc|rr|r}
        \toprule
        Problem & $\lambda$ & $\alpha$ & Scenarios & \textbf{QNN} & \textbf{IQNN} & \textbf{SAA} \\
        \midrule
        \multirow{12}{*}{CFLP-10-10} & \multirow{4}{*}{0.1} & \multirow{2}{*}{0.7}  & 500 & 8,329.23 (6.26\%) & 7,927.49 (1.14\%) & \textbf{7,838.46} \\ \multirow{12}{*}{(minimize)}
                                     &                      &                       & 1000 & 8,298.88 (6.62\%) & 7,901.39 (1.51\%) & \textbf{7,783.66} \\
                                     \cmidrule(lr){3-7}
                                     &                      & \multirow{2}{*}{0.9}  & 500 & 8,492.40 (5.82\%) & 8,114.51 (1.11\%) & \textbf{8,025.25} \\
                                     &                      &                       & 1000 & 8,453.05 (5.86\%) & 8,083.66 (1.23\%) & \textbf{7,985.19} \\
                                     \cmidrule(lr){2-7}
                                     & \multirow{4}{*}{0.5} & \multirow{2}{*}{0.7}  & 500 & 11,005.07 (0.33\%) & 11,177.91 (1.91\%) & \textbf{10,968.72} \\
                                     &                      &                       & 1000 & 10,984.90 (0.77\%) & 11,126.14 (2.07\%) & \textbf{10,900.81} \\
                                     \cmidrule(lr){3-7}
                                     &                      & \multirow{2}{*}{0.9}  & 500 & 11,924.57 (2.49\%) & 12,113.76 (4.12\%) & \textbf{11,634.79} \\
                                     &                      &                       & 1000 & 11,944.44 (2.87\%) & 12,037.83 (3.67\%) & \textbf{11,611.15} \\
                                     \cmidrule(lr){2-7}
                                     & \multirow{4}{*}{1} & \multirow{2}{*}{0.7}  & 500 & 14,961.60 (1.10\%) & 15,241.04 (2.99\%) & \textbf{14,798.28} \\
                                     &                      &                       & 1000 & 14,928.63 (1.02\%) & 15,157.14 (2.57\%) & \textbf{14,778.02} \\
                                     \cmidrule(lr){3-7}
                                     &                      & \multirow{2}{*}{0.9}  & 500 & 16,246.62 (0.72\%) & \textbf{16,112.80} (-0.11\%) & 16,130.39 \\
                                     &                      &                       & 1000 & 16,121.64 (0.21\%) & \textbf{16,040.33} (-0.30\%) & 16,088.48 \\
        \midrule
        \multirow{12}{*}{CFLP-25-25} & \multirow{4}{*}{0.1} & \multirow{2}{*}{0.7}  & 500 & 13,401.57 (-0.18\%) & \textbf{13,169.96} (-1.90\%) & 13,425.15 \\ \multirow{12}{*}{(minimize)} 
                                     &                      &                       & 1000 & 13,376.71 (0.70\%) & \textbf{13,074.58} (-1.57\%) & 13,283.26 \\
                                     \cmidrule(lr){3-7}
                                     &                      & \multirow{2}{*}{0.9}  & 500 & 13,545.35 (0.42\%) & \textbf{13,325.99} (-1.21\%) & 13,488.87 \\
                                     &                      &                       & 1000 & 13,544.26 (-4.00\%) & \textbf{13,288.85} (-5.81\%) & 14,108.71 \\
                                     \cmidrule(lr){2-7}
                                     & \multirow{4}{*}{0.5} & \multirow{2}{*}{0.7}  & 500 & 19,123.39 (3.11\%) & 18,630.72 (0.45\%) & \textbf{18,546.91} \\
                                     &                      &                       & 1000 & 19,074.51 (-1.89\%) & \textbf{18,438.73} (-5.16\%) & 19,442.64 \\
                                     \cmidrule(lr){3-7}
                                     &                      & \multirow{2}{*}{0.9}  & 500 & 20,034.65 (2.90\%) & 19,537.58 (0.35\%) & \textbf{19,470.00} \\
                                     &                      &                       & 1000 & 20,037.17 (1.46\%) & \textbf{19,507.29} (-1.22\%) & 19,748.53 \\
                                     \cmidrule(lr){2-7}
                                     & \multirow{4}{*}{1} & \multirow{2}{*}{0.7}  & 500 & 26,275.65 (6.39\%) & 25,456.68 (3.08\%) & \textbf{24,696.49} \\
                                     &                      &                       & 1000 & 26,196.76 (-0.27\%) & \textbf{25,143.93} (-4.28\%) & 26,266.91 \\
                                     \cmidrule(lr){3-7}
                                     &                      & \multirow{2}{*}{0.9}  & 500 & 27,496.57 (5.03\%) & 27,302.06 (4.28\%) & \textbf{26,180.85} \\
                                     &                      &                       & 1000 & 27,448.71 (4.26\%) & 27,280.35 (3.62\%) & \textbf{26,326.05} \\
        \midrule
        \multirow{12}{*}{CFLP-50-50} & \multirow{4}{*}{0.1} & \multirow{2}{*}{0.7}  & 500 & \textbf{29,060.30} (-2.60\%) & 30,604.24 (2.58\%) & 29,834.87 \\ \multirow{12}{*}{(minimize)}
                                     &                      &                       & 1000 & \textbf{28,900.38} (-14.32\%) & 30,442.68 (-9.75\%) & 33,729.72 \\
                                     \cmidrule(lr){3-7}
                                     &                      & \multirow{2}{*}{0.9}  & 500 & \textbf{29,367.16} (-4.36\%) & 30,905.23 (0.65\%) & 30,706.34 \\
                                     &                      &                       & 1000 & \textbf{29,210.45} (-3.58\%) & 30,735.29 (1.45\%) & 30,295.92 \\
                                     \cmidrule(lr){2-7}
                                     & \multirow{4}{*}{0.5} & \multirow{2}{*}{0.7}  & 500 & \textbf{41,483.83} (-1.90\%) & 43,724.37 (3.39\%) & 42,288.94 \\
                                     &                      &                       & 1000 & \textbf{41,137.64} (-9.40\%) & 43,499.44 (-4.20\%) & 45,405.05 \\
                                     \cmidrule(lr){3-7}
                                     &                      & \multirow{2}{*}{0.9}  & 500 & \textbf{43,071.39} (-2.63\%) & 45,229.29 (2.25\%) & 44,232.52 \\
                                     &                      &                       & 1000 & \textbf{42,673.89} (-9.77\%) & 44,962.50 (-4.93\%) & 47,294.99 \\
                                     \cmidrule(lr){2-7}
                                     & \multirow{4}{*}{1} & \multirow{2}{*}{0.7}    & 500 & 56,223.62 (5.40\%) & 60,124.52 (12.71\%) & \textbf{53,345.19} \\
                                     &                      &                       & 1000 & \textbf{55,934.78} (-3.28\%) & 59,820.39 (3.44\%) & 57,829.85 \\
                                     \cmidrule(lr){3-7}
                                     &                      & \multirow{2}{*}{0.9}  & 500 & \textbf{59,242.84} (-1.28\%) & 63,134.36 (5.21\%) & 60,010.53 \\
                                     &                      &                       & 1000 & \textbf{58,806.49} (-2.00\%) & 62,746.51 (4.56\%) & 60,009.27 \\
        \midrule
        \multirow{12}{*}{IP-I-H}     & \multirow{4}{*}{0.1} & \multirow{2}{*}{0.7}  & 1461 & \textbf{70.41} (1.51\%) & \textbf{70.41} (1.51\%) & 69.39 \\\multirow{12}{*}{(maximize)}
                                     &                      &                       & 10000 & \textbf{70.33} (9.21\%) & \textbf{70.33} (9.21\%) & 64.40 \\
                                     \cmidrule(lr){3-7}
                                     &                      & \multirow{2}{*}{0.9}  & 1461 & \textbf{69.56} (1.19\%) & \textbf{69.56} (1.19\%) & 68.74 \\
                                     &                      &                       & 10000 & \textbf{69.85} (9.07\%) & \textbf{69.85} (9.07\%) & 64.04 \\
                                     \cmidrule(lr){2-7}
                                     & \multirow{4}{*}{0.5} & \multirow{2}{*}{0.7}  & 1461 & 87.77 (-0.86\%) & 87.77 (-0.86\%) & \textbf{88.53} \\
                                     &                      &                       & 10000 & \textbf{88.19} (10.35\%) & \textbf{88.19} (10.35\%) & 79.92 \\
                                     \cmidrule(lr){3-7}
                                     &                      & \multirow{2}{*}{0.9}  & 1461 & 85.16 (-0.32\%) & 85.16 (-0.32\%) & \textbf{85.43} \\
                                     &                      &                       & 10000 & \textbf{85.31} (3.24\%) & \textbf{85.31} (3.24\%) & 82.63 \\
                                     \cmidrule(lr){2-7}
                                     & \multirow{4}{*}{1} & \multirow{2}{*}{0.7}    & 1461 & \textbf{110.18} (0.89\%) & \textbf{110.18} (0.89\%) & 109.21 \\
                                     &                      &                       & 10000 & \textbf{110.59} (8.18\%) & \textbf{110.59} (8.18\%) & 102.23 \\
                                     \cmidrule(lr){3-7}
                                     &                      & \multirow{2}{*}{0.9}  & 1461 & \textbf{104.58} (2.46\%) & \textbf{104.58} (2.46\%) & 102.07 \\
                                     &                      &                       & 10000 & \textbf{104.80} (11.49\%) & \textbf{104.80} (11.49\%) & 94.00 \\
        \bottomrule
    \end{tabular}
    \end{adjustbox}
    \caption{True objective results for risk-averse optimization. Relative differences between the QNN and IQNN objective values and the SAA approach are shown in brackets.}
    \label{tab:true_obj_averse}
\end{table}

We similarly observe promising results for the quantile-based framework in its ability to produce high-quality solutions in risk-averse optimization. In the worst cases, QNN and IQNN present a gap of 6.62\% and 12.71\% to the best solution found using SAA. On the other hand, in the best cases, the objective is improved by 14.32\% and 11.49\%. The only optimization problem in which SAA performs generally better is the CFLP-10-10. This may be due to the fact that the problem is smaller and so is easier to solve the problem in an extensive form. However, in the rest of the two-stage problems, either QNN or IQNN is generally able to obtain better solutions than the SAA formulation. This provides a great insight into the usefulness of the developed framework, especially when we are trying to solve big problems, with a considerably large number of scenarios, or for different risk-aversion levels.

\subsection{Quantile crossing and prescriptive impact}

We conclude the experimental section by analyzing the prescriptive impact of the quantile crossing tolerance parameter $\Delta$ for the QNN surrogate methodology (see Section \ref{sec:qnn_prob_formulation}). 

Table \ref{tab:delta_dif} provides a comprehensive comparison of key metrics when solving three different two-stage stochastic Capacitated Facility Location Problems (CFLP) under different risk-aversion levels. We give the true objective of the solution (the heuristic solution evaluated in a specific scenario set), the solving times in seconds, and the number of nodes explored by the solver for the surrogate formulation. The true objective value has been scaled to reveal the relative gap between solutions with different tolerance parameter values. Modifying the value of $\Delta$ reveals the impact of the parameter in the modeling and mixed-integer optimization process.

\begin{table}
    \centering
    \begin{adjustbox}{width=\textwidth} 
    \begin{tabular}{cc|ccccc|c}
        \toprule
        Problem &  Metric & $\Delta = 0$ & $\Delta = 10$ & $\Delta = 50$ & $\Delta = 100$ & $\Delta = 500$ & No Constraint \\
        \midrule
        \multirow{3}{*}{\thead{CFLP-10-10\\($\lambda=0$)}} & True Obj. & 1.01 & 1.00 & 1.00 & 1.00 & 1.00 & 1.00 \\
                                                    & Solving Time & 2.16 & 1.70 & 1.73 & 1.72 & 2.01 & 0.57 \\
                                                    & Nodes Explored & 4,954 & 1,053 & 2,501 & 2,117 & 3,262 & 682 \\
        \midrule
        \multirow{3}{*}{\thead{CFLP-25-25\\($\lambda=0.1, \sigma=0.9$)}} & True Obj. & 1.17 & 1.16 & 1.02 & 1.00 & 1.17 & 1.18 \\
                                                    & Solving Time & 0.40 & 0.30 & 0.22 & 0.24 & 0.20 & 0.02 \\
                                                    & Nodes Explored & 2,183 & 1,615 & 121 & 22 & 469 & 1 \\
        \midrule
        \multirow{3}{*}{\thead{CFLP-50-50\\($\lambda=0.5, \sigma=0.7$)}} & True Obj. & 1.02 & 1.00 & 1.01 & 1.02 & 1.03 & 1.03 \\
                                                    & Solving Time & 3.61 & 3.12 & 2.92 & 2.58 & 15.88 & 7.39 \\
                                                    & Nodes Explored & 34,377 & 30,428 & 26,960 & 21,360 & 53,865 & 106,224 \\
        \bottomrule
    \end{tabular}
    \end{adjustbox}
    \caption{Differences in the true objective, solving times, and nodes explored when employing a QNN surrogate model for different levels of tolerance parameter $\Delta$.}
    \label{tab:delta_dif}
\end{table}

Regarding the relative value of the true objective, we can see that in general, this remains close to 1, showing the effectiveness of the QNN surrogate model regardless of the tolerance $\Delta$ that we set. Only for the CFLP-25-25 problem do we observe a slight increase in the true objective when $\Delta$ is set to be either too small or large, suggesting a small degradation in solution quality (but no more than 18\% for this problem).

For the solving times and the number of nodes explored, there is also no clear pattern as to whether increasing or decreasing the value of the tolerance parameter $\Delta$ helps in faster solution of the problem using a MILP solver. However, it seems clear that completely omitting the quantile crossing constraint may lead to a higher number of nodes explored in larger problems, such as CFLP-50-50, and more difficulties in finding the optimal solution owing to the larger feasible space.

In short, adding the quantile crossing constraint is generally beneficial for the solving time and solution quality (especially in larger optimization problems). The tuning of the $\Delta$ parameter can lead to better solutions in practice, so it may be worth spending some time on this task. As we showed in the previous sections, this parameter selection is relatively fast to complete.

\section{Conclusions}
\label{sec:conclu}

Two-stage stochastic programs are typically solved using the sample average approximation (SAA) approach. 
Nevertheless, the curse of dimensionality in second-stage variables and constraints present in SAA formulations represents a key challenge for the optimization community. Data-driven and decomposition techniques have been proven useful to speed up the solving of these problems, but they rely on problem-specific solution algorithms.

This paper introduces an innovative Quantile Neural Network-based framework as an alternative approximation for two-stage stochastic programs. We show that the proposed framework is fast for data generation, training of the quantile neural network, and solving the surrogate problem with the neural network embedded as a second-stage approximation. Moreover, the incorporation of the quantile distribution enhances the framework's versatility by facilitating the consideration of risk measures, such as Conditional Value at Risk (CVaR), in addition to the standard expected value approximations. Two neural network architectures are introduced, a classical feed-forward quantile network (QNN) and an incremental network that ensures non-crossing quantiles (IQNN).

An extensive computational case study is carried out for both risk-neutral and risk-averse optimization problems. The framework exhibits a high versatility, obtaining good results across several problems considered, across different risk-aversion levels, and across different scenario-set sizes. The framework especially stands out for big problems with a high number of first and second-stage decision variables, where the SAA quickly becomes intractable. In these kinds of problems, a solution can be obtained in less than one second using a pre-trained (I)QNN, and is often a better solution than one using the SAA approach, even after multiple hours of running time. 

Future research could include the improvement of cost tail modeling in problems with a high risk-aversion level. This could be made, e.g., with a double network embedding for both the complete distribution that approximates the expected value, and a tail-focused network that approximates the CVaR. Another study could be made regarding the online adjustment of the surrogate model when we have more precise information on the stochastic scenarios or covariates for the problem.

\section*{Credit authorship contribution statement}

\textbf{Antonio Alc\'antara:} Conceptualization, Data Curation, Investigation, Methodology, Software, Writing - Original Draft. \textbf{Carlos Ruiz:} Conceptualization, Funding acquisition, Supervision, Writing - Original Draft, Writing - Review \& Editing. \textbf{Calvin Tsay:} Conceptualization, Funding acquisition, Methodology, Supervision, Writing - Original Draft, Writing - Review \& Editing.

\section*{Declaration of competing interest}

The authors declare that they have no known competing financial interests or personal relationships that could have appeared to influence the work reported in this paper.

\section*{Acknowledgements}
The authors gratefully acknowledge the financial support from MCIN/AEI/ 10.13039/ 501100011033, project PID2020-116694GB-I00, and the FPU grant (FPU20/00916).


\bibliographystyle{model5-names}\biboptions{authoryear}
\bibliography{biblio}

\newpage
\appendix
\section{Model Selection}
\label{appendix:model_sel}

\subsection{Hyper-parameter search space}

Hyper-parameters and model architectures for both QNN and IQNN have been selected with a random search of 100 potential configurations, summarized in Table \ref{tab:random_search_space}. The numbers of hidden layers and epochs were fixed to one and 2,000, respectively. The optimizer, the number of neurons, and batch size were selected across a pool of values using grid search. The learning rate and dropout proportion were sampled from a uniform distribution with the given bounds.

\begin{table}[ht]
    \centering
    \begin{tabular}{l|c}
        \toprule
                 \textbf{Hyper-parameter} & \textbf{(I)QNN} \\
        \midrule
        Batch size & $\{64, 128, 256, 512\}$  \\
        Learning rate & $\left[1e^{-5}, 1e^{-1} \right]$ \\
        Optimizer & $\{\text{Adam}, \text{Adagrad}, \text{RMSprop}\}$ \\
        Dropout & $\left[0, 0.30 \right]$  \\
        \# Epochs & 2,000  \\
        \# Neurons per hidden layer & $\{32, 64, 128, 256\}$  \\
        \bottomrule
    \end{tabular}
    \caption{Hyper-parameter search space for both QNN and IQNN models.}
    \label{tab:random_search_space}
\end{table}

\subsection{(I)QNN selected models}

The selection of the best hyper-parameter configuration is based on minimizing the validation quantile loss. The validation set is built with the 20\% of the complete dataset (20,000 samples). The best configurations for both QNN and IQNN models in each optimization problem setting are shown in Table \ref{tab:best_conf}. Note that these configurations remain the same regardless of the risk-aversion level or the quantile of interest in the optimization problem.

\begin{table}[ht]
    \centering
    \begin{adjustbox}{width=1\textwidth}
    \begin{tabular}{ll|cccc}
        \toprule
        \textbf{Model} & \textbf{Hyper-parameter} & CFLP-10-10 & CFLP-25-25 & CFLP-50-50 & IP-I-H \\
        \midrule
        \multirow{5}{*}{QNN} & Batch size & $64$ & $128$ & $256$ & $256$ \\
                             & Learning rate & $0.0358$ & $0.0088$ & $0.0008$ & $0.0037$\\
                             & Optimizer & Adagrad & Adam & Adam & RMSprop \\
                             & Dropout & $0.0979$ & $0.0371$ & $0.0320$ & $0.0025$ \\
                             & \# Neurons per hidden layer & $256$ & $32$ & $32$ & $32$\\
        \midrule
        \multirow{5}{*}{IQNN} & Batch size & $128$ & $256$ & $512$ & $512$ \\
                             & Learning rate & $0.0093$ & $0.0001$ & $0.0001$ & $0.0014$ \\
                             & Optimizer & Adam & Adam & Adam & Adam \\
                             & Dropout & $0.0479$ & $0.0003$ & $0.0022$ & $0.0000$\\
                             & \# Neurons per hidden layer & $64$ & $256$ & $32$ & $128$ \\
        \bottomrule
    \end{tabular}
    \end{adjustbox}
    \caption{QNN and IQNN best hyper-parameter configurations.}
    \label{tab:best_conf}
\end{table}

\subsection{Quantile-crossing tolerance selection}

As stated throughout the article, the selection of the quantile-crossing tolerance parameter $\Delta$ for the QNN surrogate approach was done in a prescriptive way (see Algorithm \ref{alg:tolerance_selection}). For this selection, we evaluated the values of 0, 10, 50, 100, and 500 for $\Delta$. In addition, we studied the solution when the quantile-crossing constraint was not added to the surrogate problem (\ref{eq:quant2sp_cons8}), though we found this to never be the best option to solve the surrogate problem.

Following Algorithm \ref{alg:tolerance_selection}, a scenario set of size 50 was employed for the CFLP-10-10, CFLP-25-25, and IP-I-H optimization problems, while a size of 30 was used for CFLP-50-50. The optimal value of $\Delta$ ranged between $0$ and $50$ for CFLP-10-10, $10$ and $100$ for CFLP-25-25 and CFLP-50-50, and $10$ and $50$ for IP-I-H problem.

\end{document}